\documentclass[a4paper,11pt]{amsart}
\usepackage{amsmath,mathtools,amscd,amssymb}
\usepackage[T1]{fontenc}
\usepackage[utf8]{inputenc}
\usepackage[english]{babel}
\usepackage{graphicx}

\begin{document}
\newtheorem{Lem}{Lemma \thesection.}
\newtheorem{Th}[Lem]{Theorem \thesection.}
\newtheorem{Cor}[Lem]{Corollary \thesection.}
\newtheorem{Def}[Lem]{Definition \thesection.}
\newtheorem{Ex}[Lem]{Examples \thesection.}
\newtheorem{Prop}[Lem]{Proposition \thesection.}
\newtheorem{Rem}[Lem]{Remark \thesection.}
\newtheorem{Lem and Def}[Lem]{Lemma and Definition \thesection.}
\setcounter{Lem}{0}
\def\cal{\mathcal}
\def\bb{\mathbb} 
\def\a{\alpha }
\def\b{\beta }
\def\g{\gamma }
\def\d{\delta }
\def\D{\Delta }
\def\e{\epsilon }
\def\f{\varphi }
\def\F{\Phi }
\def\g{\gamma }
\def\G{\Gamma }
\def\k{\kappa }
\def\l{\lambda }
\def\m{\mu }
\def\o{\omega }
\def\O{\Omega }
\def\p{\pi }
\def\P{\Pi }
\def\s{\sigma }
\def\S{\Sigma }
\def\t{\theta }
\def\T{\Theta }
\def\z{\zeta }
\def\dim{\rm dim\; }
\def\codim{\rm codim\; }
\def\iff{ if and only if }
\def\ker{{\rm Ker\,}}
\def\im{{\rm Im\,}}
\def\coker{{\rm Coker\,}}
\def\Dem{D\'emonstration: }
\def\Card{{\rm Card\ }}
\def\ot{\otimes }
\def\part{\partial }
\def\tr{\rm tr}
\def\dps{\displaystyle }




\newcommand{\tlowername}[2]%
{$\stackrel{\makebox[1pt]{#1}}%
{\begin{picture}(0,0)%
\put(0,0){\makebox(0,6)[t]{\makebox[1pt]{$#2$}}}%
\end{picture}}$}%

\newcommand{\tcase}[1]{\makebox[23pt]%
{\raisebox{2.5pt}{#1{20}}}}%

\newcommand{\Tcase}[2]{\makebox[23pt]%
{\raisebox{2.5pt}{$\stackrel{#2}{#1{20}}$}}}%

\newcommand{\tbicase}[1]{\makebox[23pt]%
{\raisebox{1pt}{#1{20}}}}%

\newcommand{\Tbicase}[3]{\makebox[23pt]{\raisebox{-7pt}%
{$\stackrel{#2}{\mbox{\tlowername{#1{20}}{\scriptstyle{#3}}}}$}}}%


\newcommand{\AR}[1]%
{\begin{picture}(#1,0)%
\put(0,0){\vector(1,0){#1}}%
\end{picture}}%

\newcommand{\DOTAR}[1]%
{\NUMBEROFDOTS=#1%
\divide\NUMBEROFDOTS by 3%
\begin{picture}(#1,0)%
\multiput(0,0)(3,0){\NUMBEROFDOTS}{\circle*{1}}%
\put(#1,0){\vector(1,0){0}}%
\end{picture}}%

\newcommand{\MONO}[1]%
{\begin{picture}(#1,0)%
\put(0,0){\vector(1,0){#1}}%
\put(2,-2){\line(0,1){4}}%
\end{picture}}%

\newcommand{\EPI}[1]%
{\begin{picture}(#1,0)(-#1,0)%
\put(-#1,0){\vector(1,0){#1}}%
\put(-6,-2){\line(0,1){4}}%
\end{picture}}%

\newcommand{\BIMO}[1]%
{\begin{picture}(#1,0)(-#1,0)%
\put(-#1,0){\vector(1,0){#1}}%
\put(-6,-2){\line(0,1){4}}%
\put(-#1,-2){\hspace{2pt}\line(0,1){4}}%
\end{picture}}%

\newcommand{\BIAR}[1]%
{\begin{picture}(#1,4)%
\put(0,0){\vector(1,0){#1}}%
\put(0,4){\vector(1,0){#1}}%
\end{picture}}%

\newcommand{\EQL}[1]%
{\begin{picture}(#1,0)%
\put(0,1){\line(1,0){#1}}%
\put(0,-1){\line(1,0){#1}}%
\end{picture}}%

\newcommand{\ADJAR}[1]%
{\begin{picture}(#1,4)%
\put(0,0){\vector(1,0){#1}}%
\put(#1,4){\vector(-1,0){#1}}%
\end{picture}}%


\newcommand{\ar}{\tcase{\AR}}%

\newcommand{\Ar}[1]{\Tcase{\AR}{#1}}%

\newcommand{\dotar}{\tcase{\DOTAR}}%

\newcommand{\Dotar}[1]{\Tcase{\DOTAR}{#1}}%

\newcommand{\mono}{\tcase{\MONO}}%

\newcommand{\Mono}[1]{\Tcase{\MONO}{#1}}%

\newcommand{\epi}{\tcase{\EPI}}%

\newcommand{\Epi}[1]{\Tcase{\EPI}{#1}}%

\newcommand{\bimo}{\tcase{\BIMO}}%

\newcommand{\Bimo}[1]{\Tcase{\BIMO}{#1}}%

\newcommand{\iso}{\Tcase{\AR}{\cong}}%

\newcommand{\Iso}[1]{\Tcase{\AR}{\cong{#1}}}%

\newcommand{\biar}{\tbicase{\BIAR}}%

\newcommand{\Biar}[2]{\Tbicase{\BIAR}{#1}{#2}}%

\newcommand{\eql}{\tcase{\EQL}}%

\newcommand{\adjar}{\tbicase{\ADJAR}}%

\newcommand{\Adjar}[2]{\Tbicase{\ADJAR}{#1}{#2}}%


\newcommand{\BKAR}[1]%
{\begin{picture}(#1,0)%
\put(#1,0){\vector(-1,0){#1}}%
\end{picture}}%

\newcommand{\BKDOTAR}[1]%
{\NUMBEROFDOTS=#1%
\divide\NUMBEROFDOTS by 3%
\begin{picture}(#1,0)%
\multiput(#1,0)(-3,0){\NUMBEROFDOTS}{\circle*{1}}%
\put(0,0){\vector(-1,0){0}}%
\end{picture}}%

\newcommand{\BKMONO}[1]%
{\begin{picture}(#1,0)(-#1,0)%
\put(0,0){\vector(-1,0){#1}}%
\put(-2,-2){\line(0,1){4}}%
\end{picture}}%

\newcommand{\BKEPI}[1]%
{\begin{picture}(#1,0)%
\put(#1,0){\vector(-1,0){#1}}%
\put(6,-2){\line(0,1){4}}%
\end{picture}}%

\newcommand{\BKBIMO}[1]%
{\begin{picture}(#1,0)%
\put(#1,0){\vector(-1,0){#1}}%
\put(6,-2){\line(0,1){4}}%
\put(#1,-2){\hspace{-2pt}\line(0,1){4}}%
\end{picture}}%

\newcommand{\BKBIAR}[1]%
{\begin{picture}(#1,4)%
\put(#1,0){\vector(-1,0){#1}}%
\put(#1,4){\vector(-1,0){#1}}%
\end{picture}}%

\newcommand{\BKADJAR}[1]%
{\begin{picture}(#1,4)%
\put(0,4){\vector(1,0){#1}}%
\put(#1,0){\vector(-1,0){#1}}%
\end{picture}}%


\newcommand{\bkar}{\tcase{\BKAR}}%

\newcommand{\Bkar}[1]{\Tcase{\BKAR}{#1}}%

\newcommand{\bkdotar}{\tcase{\BKDOTAR}}%

\newcommand{\Bkdotar}[1]{\Tcase{\BKDOTAR}{#1}}%

\newcommand{\bkmono}{\tcase{\BKMONO}}%

\newcommand{\Bkmono}[1]{\Tcase{\BKMONO}{#1}}%

\newcommand{\bkepi}{\tcase{\BKEPI}}%

\newcommand{\Bkepi}[1]{\Tcase{\BKEPI}{#1}}%

\newcommand{\bkbimo}{\tcase{\BKBIMO}}%

\newcommand{\Bkbimo}[1]{\Tcase{\BKBIMO}{\hspace{9pt}#1}}%

\newcommand{\bkiso}{\Tcase{\BKAR}{\cong}}%

\newcommand{\Bkiso}[1]{\Tcase{\BKAR}{\cong{#1}}}%

\newcommand{\bkbiar}{\tbicase{\BKBIAR}}%

\newcommand{\Bkbiar}[2]{\Tbicase{\BKBIAR}{#1}{#2}}%

\newcommand{\bkeql}{\tcase{\EQL}}%

\newcommand{\bkadjar}{\tbicase{\BKADJAR}}%

\newcommand{\Bkadjar}[2]{\Tbicase{\BKADJAR}{#1}{#2}}%


\newcommand{\lowername}[2]%
{$\stackrel{\makebox[1pt]{#1}}%
{\begin{picture}(0,0)%
\truex{600}%
\put(0,0){\makebox(0,\value{x})[t]{\makebox[1pt]{$#2$}}}%
\end{picture}}$}%

\newcommand{\hcase}[2]%
{\makebox[0pt]%
{\raisebox{-1pt}[0pt][0pt]{#1{#2}}}}%

\newcommand{\Hcase}[3]%
{\makebox[0pt]
{\raisebox{-1pt}[0pt][0pt]%
{$\stackrel{\makebox[0pt]{$\textstyle{#2}$}}{#1{#3}}$}}}%

\newcommand{\hcasE}[3]%
{\makebox[0pt]%
{\raisebox{-9pt}[0pt][0pt]%
{\lowername{#1{#3}}{#2}}}}%

\newcommand{\hbicase}[2]%
{\makebox[0pt]%
{\raisebox{-2.5pt}[0pt][0pt]{#1{#2}}}}%

\newcommand{\Hbicase}[4]%
{\makebox[0pt]
{\raisebox{-10.5pt}[0pt][0pt]%
{$\stackrel{\makebox[0pt]{$\textstyle{#2}$}}%
{\mbox{\lowername{#1{#4}}{#3}}}$}}}%


\newcommand{\EAR}[1]%
{\begin{picture}(#1,0)%
\put(0,0){\vector(1,0){#1}}%
\end{picture}}%

\newcommand{\EDOTAR}[1]%
{\truex{100}\truey{300}%
\NUMBEROFDOTS=#1%
\divide\NUMBEROFDOTS by \value{y}%
\begin{picture}(#1,0)%
\multiput(0,0)(\value{y},0){\NUMBEROFDOTS}%
{\circle*{\value{x}}}%
\put(#1,0){\vector(1,0){0}}%
\end{picture}}%

\newcommand{\EMONO}[1]%
{\begin{picture}(#1,0)%
\put(0,0){\vector(1,0){#1}}%
\truex{300}\truey{600}%
\put(\value{x},-\value{x}){\line(0,1){\value{y}}}%
\end{picture}}%

\newcommand{\EEPI}[1]%
{\begin{picture}(#1,0)(-#1,0)%
\put(-#1,0){\vector(1,0){#1}}%
\truex{300}\truey{600}\truez{800}%
\put(-\value{z},-\value{x}){\line(0,1){\value{y}}}%
\end{picture}}%

\newcommand{\EBIMO}[1]%
{\begin{picture}(#1,0)(-#1,0)%
\put(-#1,0){\vector(1,0){#1}}%
\truex{300}\truey{600}\truez{800}%
\put(-\value{z},-\value{x}){\line(0,1){\value{y}}}%
\put(-#1,-\value{x}){\hspace{3pt}\line(0,1){\value{y}}}%
\end{picture}}%

\newcommand{\EBIAR}[1]%
{\truex{400}%
\begin{picture}(#1,\value{x})%
\put(0,0){\vector(1,0){#1}}%
\put(0,\value{x}){\vector(1,0){#1}}%
\end{picture}}%

\newcommand{\EEQL}[1]%
{\begin{picture}(#1,0)%
\truex{200}%
\put(0,\value{x}){\line(1,0){#1}}%
\put(0,0){\line(1,0){#1}}%
\end{picture}}%

\newcommand{\EADJAR}[1]%
{\truex{400}%
\begin{picture}(#1,\value{x})%
\put(0,0){\vector(1,0){#1}}%
\put(#1,\value{x}){\vector(-1,0){#1}}%
\end{picture}}%


\newcommand{\earv}[1]{\hcase{\EAR}{#100}}%

\newcommand{\ear}%
{\hspace{\SOURCE\unitlength}%
\hcase{\EAR}{\ARROWLENGTH}}%

\newcommand{\Earv}[2]{\Hcase{\EAR}{#1}{#200}}%

\newcommand{\Ear}[1]%
{\hspace{\SOURCE\unitlength}%
\Hcase{\EAR}{#1}{\ARROWLENGTH}}%

\newcommand{\eaRv}[2]{\hcasE{\EAR}{#1}{#200}}%

\newcommand{\eaR}[1]%
{\hspace{\SOURCE\unitlength}%
\hcasE{\EAR}{#1}{\ARROWLENGTH}}%

\newcommand{\edotarv}[1]{\hcase{\EDOTAR}{#100}}%

\newcommand{\edotar}%
{\hspace{\SOURCE\unitlength}%
\hcase{\EDOTAR}{\ARROWLENGTH}}%

\newcommand{\Edotarv}[2]{\Hcase{\EDOTAR}{#1}{#200}}%

\newcommand{\Edotar}[1]%
{\hspace{\SOURCE\unitlength}%
\Hcase{\EDOTAR}{#1}{\ARROWLENGTH}}%

\newcommand{\edotaRv}[2]{\hcasE{\EDOTAR}{#1}{#200}}%

\newcommand{\edotaR}[1]%
{\hspace{\SOURCE\unitlength}%
\hcasE{\EDOTAR}{#1}{\ARROWLENGTH}}%

\newcommand{\emonov}[1]{\hcase{\EMONO}{#100}}%

\newcommand{\emono}%
{\hspace{\SOURCE\unitlength}%
\hcase{\EMONO}{\ARROWLENGTH}}%

\newcommand{\Emonov}[2]{\Hcase{\EMONO}{#1}{#200}}%

\newcommand{\Emono}[1]%
{\hspace{\SOURCE\unitlength}%
\Hcase{\EMONO}{#1}{\ARROWLENGTH}}%

\newcommand{\emonOv}[2]{\hcasE{\EMONO}{#1}{#200}}%

\newcommand{\emonO}[1]%
{\hspace{\SOURCE\unitlength}%
\hcasE{\EMONO}{#1}{\ARROWLENGTH}}%

\newcommand{\eepiv}[1]{\hcase{\EEPI}{#100}}%

\newcommand{\eepi}%
{\hspace{\SOURCE\unitlength}%
\hcase{\EEPI}{\ARROWLENGTH}}%

\newcommand{\Eepiv}[2]{\Hcase{\EEPI}{#1}{#200}}%

\newcommand{\Eepi}[1]%
{\hspace{\SOURCE\unitlength}%
\Hcase{\EEPI}{#1}{\ARROWLENGTH}}%

\newcommand{\eepIv}[2]{\hcasE{\EEPI}{#1}{#200}}%

\newcommand{\eepI}[1]%
{\hspace{\SOURCE\unitlength}%
\hcasE{\EEPI}{#1}{\ARROWLENGTH}}%

\newcommand{\ebimov}[1]{\hcase{\EBIMO}{#100}}%

\newcommand{\ebimo}%
{\hspace{\SOURCE\unitlength}%
\hcase{\EBIMO}{\ARROWLENGTH}}%

\newcommand{\Ebimov}[2]{\Hcase{\EBIMO}{#1}{#200}}%

\newcommand{\Ebimo}[1]%
{\hspace{\SOURCE\unitlength}%
\Hcase{\EBIMO}{#1}{\ARROWLENGTH}}%

\newcommand{\ebimOv}[2]{\hcasE{\EBIMO}{#1}{#200}}%

\newcommand{\ebimO}[1]%
{\hspace{\SOURCE\unitlength}%
\hcasE{\EBIMO}{#1}{\ARROWLENGTH}}%

\newcommand{\eisov}[1]{\Hcase{\EAR}{\cong}{#100}}%

\newcommand{\eiso}%
{\hspace{\SOURCE\unitlength}%
\Hcase{\EAR}{\cong}{\ARROWLENGTH}}%

\newcommand{\Eisov}[2]{\Hcase{\EAR}{\cong#1}{#200}}%

\newcommand{\Eiso}[1]%
{\hspace{\SOURCE\unitlength}%
\Hcase{\EAR}{\cong#1}{\ARROWLENGTH}}%

\newcommand{\eisOv}[2]{\hcasE{\EAR}{\cong#1}{#200}}%

\newcommand{\eisO}[1]%
{\hspace{\SOURCE\unitlength}%
\hcasE{\EAR}{\cong#1}{\ARROWLENGTH}}%

\newcommand{\ebiarv}[1]{\hbicase{\EBIAR}{#100}}%

\newcommand{\ebiar}%
{\hspace{\SOURCE\unitlength}%
\hbicase{\EBIAR}{\ARROWLENGTH}}%

\newcommand{\Ebiarv}[3]{\Hbicase{\EBIAR}{#1}{#2}{#300}}%

\newcommand{\Ebiar}[2]%
{\hspace{\SOURCE\unitlength}%
\Hbicase{\EBIAR}{#1}{#2}{\ARROWLENGTH}}%

\newcommand{\eeqlv}[1]{\hcase{\EEQL}{#100}}%

\newcommand{\eeql}%
{\hspace{\SOURCE\unitlength}%
\hbicase{\EEQL}{\ARROWLENGTH}}%

\newcommand{\eadjarv}[1]{\hbicase{\EADJAR}{#100}}%

\newcommand{\eadjar}%
{\hspace{\SOURCE\unitlength}%
\hbicase{\EADJAR}{\ARROWLENGTH}}%

\newcommand{\Eadjarv}[3]{\Hbicase{\EADJAR}{#1}{#2}{#300}}%

\newcommand{\Eadjar}[2]%
{\hspace{\SOURCE\unitlength}%
\Hbicase{\EADJAR}{#1}{#2}{\ARROWLENGTH}}%


\newcommand{\WAR}[1]%
{\begin{picture}(#1,0)%
\put(#1,0){\vector(-1,0){#1}}%
\end{picture}}%

\newcommand{\WDOTAR}[1]%
{\truex{100}\truey{300}%
\NUMBEROFDOTS=#1%
\divide\NUMBEROFDOTS by \value{y}%
\begin{picture}(#1,0)%
\multiput(#1,0)(-\value{y},0){\NUMBEROFDOTS}%
{\circle*{\value{x}}}%
\put(0,0){\vector(-1,0){0}}%
\end{picture}}%

\newcommand{\WMONO}[1]%
{\begin{picture}(#1,0)(-#1,0)%
\put(0,0){\vector(-1,0){#1}}%
\truex{300}\truey{600}%
\put(-\value{x},-\value{x}){\line(0,1){\value{y}}}%
\end{picture}}%

\newcommand{\WEPI}[1]%
{\begin{picture}(#1,0)%
\put(#1,0){\vector(-1,0){#1}}%
\truex{300}\truey{600}\truez{800}%
\put(\value{z},-\value{x}){\line(0,1){\value{y}}}%
\end{picture}}%

\newcommand{\WBIMO}[1]%
{\begin{picture}(#1,0)%
\put(#1,0){\vector(-1,0){#1}}%
\truex{300}\truey{600}\truez{800}%
\put(\value{z},-\value{x}){\line(0,1){\value{y}}}%
\put(#1,-\value{x}){\hspace{-3pt}\line(0,1){\value{y}}}%
\end{picture}}%

\newcommand{\WBIAR}[1]%
{\truex{400}%
\begin{picture}(#1,\value{x})%
\put(#1,0){\vector(-1,0){#1}}%
\put(#1,\value{x}){\vector(-1,0){#1}}%
\end{picture}}%

\newcommand{\WADJAR}[1]%
{\truex{400}%
\begin{picture}(#1,\value{x})%
\put(0,\value{x}){\vector(1,0){#1}}%
\put(#1,0){\vector(-1,0){#1}}%
\end{picture}}%


\newcommand{\warv}[1]{\hcase{\WAR}{#100}}%

\newcommand{\war}%
{\hspace{\SOURCE\unitlength}%
\hcase{\WAR}{\ARROWLENGTH}}%

\newcommand{\Warv}[2]{\Hcase{\WAR}{#1}{#200}}%

\newcommand{\War}[1]%
{\hspace{\SOURCE\unitlength}%
\Hcase{\WAR}{#1}{\ARROWLENGTH}}%

\newcommand{\waRv}[2]{\hcasE{\WAR}{#1}{#200}}%

\newcommand{\waR}[1]%
{\hspace{\SOURCE\unitlength}%
\hcasE{\WAR}{#1}{\ARROWLENGTH}}%

\newcommand{\wdotarv}[1]{\hcase{\WDOTAR}{#100}}%

\newcommand{\wdotar}%
{\hspace{\SOURCE\unitlength}%
\hcase{\WDOTAR}{\ARROWLENGTH}}%

\newcommand{\Wdotarv}[2]{\Hcase{\WDOTAR}{#1}{#200}}%

\newcommand{\Wdotar}[1]%
{\hspace{\SOURCE\unitlength}%
\Hcase{\WDOTAR}{#1}{\ARROWLENGTH}}%

\newcommand{\wdotaRv}[2]{\hcasE{\WDOTAR}{#1}{#200}}%

\newcommand{\wdotaR}[1]%
{\hspace{\SOURCE\unitlength}%
\hcasE{\WDOTAR}{#1}{\ARROWLENGTH}}%

\newcommand{\wmonov}[1]{\hcase{\WMONO}{#100}}%

\newcommand{\wmono}%
{\hspace{\SOURCE\unitlength}%
\hcase{\WMONO}{\ARROWLENGTH}}%

\newcommand{\Wmonov}[2]{\Hcase{\WMONO}{#1}{#200}}%

\newcommand{\Wmono}[1]%
{\hspace{\SOURCE\unitlength}%
\Hcase{\WMONO}{#1}{\ARROWLENGTH}}%

\newcommand{\wmonOv}[2]{\hcasE{\WMONO}{#1}{#200}}%

\newcommand{\wmonO}[1]%
{\hspace{\SOURCE\unitlength}%
\hcasE{\WMONO}{#1}{\ARROWLENGTH}}%

\newcommand{\wepiv}[1]{\hcase{\WEPI}{#100}}%

\newcommand{\wepi}%
{\hspace{\SOURCE\unitlength}%
\hcase{\WEPI}{\ARROWLENGTH}}%

\newcommand{\Wepiv}[2]{\Hcase{\WEPI}{#1}{#200}}%

\newcommand{\Wepi}[1]%
{\hspace{\SOURCE\unitlength}%
\Hcase{\WEPI}{#1}{\ARROWLENGTH}}%

\newcommand{\wepIv}[2]{\hcasE{\WEPI}{#1}{#200}}%

\newcommand{\wepI}[1]%
{\hspace{\SOURCE\unitlength}%
\hcasE{\WEPI}{#1}{\ARROWLENGTH}}%

\newcommand{\wbimov}[1]{\hcase{\WBIMO}{#100}}%

\newcommand{\wbimo}%
{\hspace{\SOURCE\unitlength}%
\hcase{\WBIMO}{\ARROWLENGTH}}%

\newcommand{\Wbimov}[2]{\Hcase{\WBIMO}{#1}{#200}}%

\newcommand{\Wbimo}[1]%
{\hspace{\SOURCE\unitlength}%
\Hcase{\WBIMO}{#1}{\ARROWLENGTH}}%

\newcommand{\wbimOv}[2]{\hcasE{\WBIMO}{#1}{#200}}%

\newcommand{\wbimO}[1]%
{\hspace{\SOURCE\unitlength}%
\hcasE{\WBIMO}{#1}{\ARROWLENGTH}}%

\newcommand{\wisov}[1]{\Hcase{\WAR}{\cong}{#100}}%

\newcommand{\wiso}%
{\hspace{\SOURCE\unitlength}%
\Hcase{\WAR}{\cong}{\ARROWLENGTH}}%

\newcommand{\Wisov}[2]{\Hcase{\WAR}{\cong#1}{#200}}%

\newcommand{\Wiso}[1]%
{\hspace{\SOURCE\unitlength}%
\Hcase{\WAR}{#1}{\ARROWLENGTH}}%

\newcommand{\wisOv}[2]{\hcasE{\WAR}{\cong#1}{#200}}%

\newcommand{\wisO}[1]%
{\hspace{\SOURCE\unitlength}%
\hcasE{\WAR}{#1}{\ARROWLENGTH}}%

\newcommand{\wbiarv}[1]{\hbicase{\WBIAR}{#100}}%

\newcommand{\wbiar}%
{\hspace{\SOURCE\unitlength}%
\hbicase{\WBIAR}{\ARROWLENGTH}}%

\newcommand{\Wbiarv}[3]{\Hbicase{\WBIAR}{#1}{#2}{#300}}%

\newcommand{\Wbiar}[2]%
{\hspace{\SOURCE\unitlength}%
\Hbicase{\WBIAR}{#1}{#2}{\ARROWLENGTH}}%

\newcommand{\weqlv}[1]{\hbicase{\EEQL}{#100}}%

\newcommand{\weql}%
{\hspace{\SOURCE\unitlength}%
\hbicase{\EEQL}{\ARROWLENGTH}}%

\newcommand{\wadjarv}[1]{\hbicase{\WADJAR}{#100}}%

\newcommand{\wadjar}%
{\hspace{\SOURCE\unitlength}%
\hbicase{\WADJAR}{\ARROWLENGTH}}%

\newcommand{\Wadjarv}[3]{\Hbicase{\WADJAR}{#1}{#2}{#300}}%

\newcommand{\Wadjar}[2]%
{\hspace{\SOURCE\unitlength}%
\Hbicase{\WADJAR}{#1}{#2}{\ARROWLENGTH}}%


\newcommand{\vcase}[2]{#1{#2}}%

\newcommand{\Vcase}[3]{\makebox[0pt]%
{\makebox[0pt][r]{\raisebox{0pt}[0pt][0pt]{${#2}\hspace{2pt}$}}}#1{#3}}%

\newcommand{\vcasE}[3]{\makebox[0pt]%
{#1{#3}\makebox[0pt][l]{\raisebox{0pt}[0pt][0pt]{\hspace{2pt}$#2$}}}}%

\newcommand{\vbicase}[2]{\makebox[0pt]{{#1{#2}}}}%

\newcommand{\Vbicase}[4]{\makebox[0pt]%
{\makebox[0pt][r]{\raisebox{0pt}[0pt][0pt]{$#2$\hspace{4pt}}}#1{#4}%
\makebox[0pt][l]{\raisebox{0pt}[0pt][0pt]{\hspace{5pt}$#3$}}}}%


\newcommand{\SAR}[1]%
{\begin{picture}(0,0)%
\put(0,0){\makebox(0,0)%
{\begin{picture}(0,#1)%
\put(0,#1){\vector(0,-1){#1}}%
\end{picture}}}\end{picture}}%

\newcommand{\SDOTAR}[1]%
{\truex{100}\truey{300}%
\NUMBEROFDOTS=#1%
\divide\NUMBEROFDOTS by \value{y}%
\begin{picture}(0,0)%
\put(0,0){\makebox(0,0)%
{\begin{picture}(0,#1)%
\multiput(0,#1)(0,-\value{y}){\NUMBEROFDOTS}%
{\circle*{\value{x}}}%
\put(0,0){\vector(0,-1){0}}%
\end{picture}}}\end{picture}}%

\newcommand{\SMONO}[1]%
{\begin{picture}(0,0)%
\put(0,0){\makebox(0,0)%
{\begin{picture}(0,#1)%
\put(0,#1){\vector(0,-1){#1}}%
\truex{300}\truey{600}%
\put(0,#1){\begin{picture}(0,0)%
\put(-\value{x},-\value{x}){\line(1,0){\value{y}}}\end{picture}}%
\end{picture}}}\end{picture}}%

\newcommand{\SEPI}[1]%
{\begin{picture}(0,0)%
\put(0,0){\makebox(0,0)%
{\begin{picture}(0,#1)%
\put(0,#1){\vector(0,-1){#1}}%
\truex{300}\truey{600}\truez{800}%
\put(-\value{x},\value{z}){\line(1,0){\value{y}}}%
\end{picture}}}\end{picture}}%

\newcommand{\SBIMO}[1]%
{\begin{picture}(0,0)%
\put(0,0){\makebox(0,0)%
{\begin{picture}(0,#1)%
\put(0,#1){\vector(0,-1){#1}}%
\truex{300}\truey{600}\truez{800}%
\put(0,#1){\begin{picture}(0,0)%
\put(-\value{x},-\value{x}){\line(1,0){\value{y}}}\end{picture}}%
\put(-\value{x},\value{z}){\line(1,0){\value{y}}}%
\end{picture}}}\end{picture}}%

\newcommand{\SBIAR}[1]%
{\begin{picture}(0,0)%
\truex{200}%
\put(0,0){\makebox(0,0)%
{\begin{picture}(0,#1)\put(-\value{x},#1){\vector(0,-1){#1}}%
\put(\value{x},#1){\vector(0,-1){#1}}%
\end{picture}}}\end{picture}}%

\newcommand{\SEQL}[1]%
{\begin{picture}(0,0)%
\truex{100}%
\put(0,0){\makebox(0,0)%
{\begin{picture}(0,#1)\put(-\value{x},#1){\line(0,-1){#1}}%
\put(\value{x},#1){\line(0,-1){#1}}%
\end{picture}}}\end{picture}}%

\newcommand{\SADJAR}[1]{\begin{picture}(0,0)%
\truex{200}%
\put(0,0){\makebox(0,0)%
{\begin{picture}(0,#1)\put(-\value{x},#1){\vector(0,-1){#1}}%
\put(\value{x},0){\vector(0,1){#1}}%
\end{picture}}}\end{picture}}%


\newcommand{\sarv}[1]{\vcase{\SAR}{#100}}%

\newcommand{\sar}{\sarv{50}}%

\newcommand{\Sarv}[2]{\Vcase{\SAR}{#1}{#200}}%

\newcommand{\Sar}[1]{\Sarv{#1}{50}}%

\newcommand{\saRv}[2]{\vcasE{\SAR}{#1}{#200}}%

\newcommand{\saR}[1]{\saRv{#1}{50}}%

\newcommand{\sdotarv}[1]{\vcase{\SDOTAR}{#100}}%

\newcommand{\sdotar}{\sdotarv{50}}%

\newcommand{\Sdotarv}[2]{\Vcase{\SDOTAR}{#1}{#200}}%

\newcommand{\Sdotar}[1]{\Sdotarv{#1}{50}}%

\newcommand{\sdotaRv}[2]{\vcasE{\SDOTAR}{#1}{#200}}%

\newcommand{\sdotaR}[1]{\sdotaRv{#1}{50}}%

\newcommand{\smonov}[1]{\vcase{\SMONO}{#100}}%

\newcommand{\smono}{\smonov{50}}%

\newcommand{\Smonov}[2]{\Vcase{\SMONO}{#1}{#200}}%

\newcommand{\Smono}[1]{\Smonov{#1}{50}}%

\newcommand{\smonOv}[2]{\vcasE{\SMONO}{#1}{#200}}%

\newcommand{\smonO}[1]{\smonOv{#1}{50}}%

\newcommand{\sepiv}[1]{\vcase{\SEPI}{#100}}%

\newcommand{\sepi}{\sepiv{50}}%

\newcommand{\Sepiv}[2]{\Vcase{\SEPI}{#1}{#200}}%

\newcommand{\Sepi}[1]{\Sepiv{#1}{50}}%

\newcommand{\sepIv}[2]{\vcasE{\SEPI}{#1}{#200}}%

\newcommand{\sepI}[1]{\sepIv{#1}{50}}%

\newcommand{\sbimov}[1]{\vcase{\SBIMO}{#100}}%

\newcommand{\sbimo}{\sbimov{50}}%

\newcommand{\Sbimov}[2]{\Vcase{\SBIMO}{#1}{#200}}%

\newcommand{\Sbimo}[1]{\Sbimov{#1}{50}}%

\newcommand{\sbimOv}[2]{\vcasE{\SBIMO}{#1}{#200}}%

\newcommand{\sbimO}[1]{\sbimOv{#1}{50}}%

\newcommand{\sisov}[1]{\vcasE{\SAR}{\cong}{#100}}%

\newcommand{\siso}{\sisov{50}}%

\newcommand{\Sisov}[2]%
{\Vbicase{\SAR}{#1\hspace{-2pt}}{\hspace{-2pt}\cong}{#200}}%

\newcommand{\Siso}[1]{\Sisov{#1}{50}}%

\newcommand{\sbiarv}[1]{\vbicase{\SBIAR}{#100}}%

\newcommand{\sbiar}{\sbiarv{50}}%

\newcommand{\Sbiarv}[3]{\Vbicase{\SBIAR}{#1}{#2}{#300}}%

\newcommand{\Sbiar}[2]{\Sbiarv{#1}{#2}{50}}%

\newcommand{\seqlv}[1]{\vbicase{\SEQL}{#100}}%

\newcommand{\seql}{\seqlv{50}}%

\newcommand{\sadjarv}[1]{\vbicase{\SADJAR}{#100}}%

\newcommand{\sadjar}{\sadjarv{50}}%

\newcommand{\Sadjarv}[3]{\Vbicase{\SADJAR}{#1}{#2}{#300}}%

\newcommand{\Sadjar}[2]{\Sadjarv{#1}{#2}{50}}%


\newcommand{\NAR}[1]%
{\begin{picture}(0,0)%
\put(0,0){\makebox(0,0)%
{\begin{picture}(0,#1)\put(0,0){\vector(0,1){#1}}%
\end{picture}}}\end{picture}}%

\newcommand{\NDOTAR}[1]%
{\truex{100}\truey{300}%
\NUMBEROFDOTS=#1%
\divide\NUMBEROFDOTS by \value{y}%
\begin{picture}(0,0)%
\put(0,0){\makebox(0,0)%
{\begin{picture}(0,#1)%
\multiput(0,0)(0,\value{y}){\NUMBEROFDOTS}%
{\circle*{\value{x}}}%
\put(0,#1){\vector(0,1){0}}%
\end{picture}}}\end{picture}}%

\newcommand{\NMONO}[1]%
{\begin{picture}(0,0)%
\put(0,0){\makebox(0,0)%
{\begin{picture}(0,#1)%
\put(0,0){\vector(0,1){#1}}%
\truex{300}\truey{600}%
\put(-\value{x},\value{x}){\line(1,0){\value{y}}}%
\end{picture}}}%
\end{picture}}%

\newcommand{\NEPI}[1]%
{\begin{picture}(0,0)%
\put(0,0){\makebox(0,0)%
{\begin{picture}(0,#1)%
\put(0,0){\vector(0,1){#1}}%
\truex{300}\truey{600}\truez{800}%
\put(0,#1){\begin{picture}(0,0)%
\put(-\value{x},-\value{z}){\line(1,0){\value{y}}}\end{picture}}%
\end{picture}}}\end{picture}}%

\newcommand{\NBIMO}[1]%
{\begin{picture}(0,0)%
\put(0,0){\makebox(0,0)%
{\begin{picture}(0,#1)%
\put(0,0){\vector(0,1){#1}}%
\truex{300}\truey{600}\truez{800}%
\put(-\value{x},\value{x}){\line(1,0){\value{y}}}%
\put(0,#1){\begin{picture}(0,0)%
\put(-\value{x},-\value{z}){\line(1,0){\value{y}}}\end{picture}}%
\end{picture}}}\end{picture}}%

\newcommand{\NBIAR}[1]%
{\begin{picture}(0,0)%
\truex{200}%
\put(0,0){\makebox(0,0)%
{\begin{picture}(0,#1)\put(-\value{x},0){\vector(0,1){#1}}%
\put(\value{x},0){\vector(0,1){#1}}%
\end{picture}}}\end{picture}}%

\newcommand{\NADJAR}[1]{\begin{picture}(0,0)%
\truex{200}%
\put(0,0){\makebox(0,0)%
{\begin{picture}(0,#1)\put(\value{x},#1){\vector(0,-1){#1}}%
\put(-\value{x},0){\vector(0,1){#1}}%
\end{picture}}}\end{picture}}%


\newcommand{\narv}[1]{\vcase{\NAR}{#100}}%

\newcommand{\nar}{\narv{50}}%

\newcommand{\Narv}[2]{\Vcase{\NAR}{#1}{#200}}%

\newcommand{\Nar}[1]{\Narv{#1}{50}}%

\newcommand{\naRv}[2]{\vcasE{\NAR}{#1}{#200}}%

\newcommand{\naR}[1]{\naRv{#1}{50}}%

\newcommand{\ndotarv}[1]{\vcase{\NDOTAR}{#100}}%

\newcommand{\ndotar}{\ndotarv{50}}%

\newcommand{\Ndotarv}[2]{\Vcase{\NDOTAR}{#1}{#200}}%

\newcommand{\Ndotar}[1]{\Ndotarv{#1}{50}}%

\newcommand{\ndotaRv}[2]{\vcasE{\NDOTAR}{#1}{#200}}%

\newcommand{\ndotaR}[1]{\ndotaRv{#1}{50}}%

\newcommand{\nmonov}[1]{\vcase{\NMONO}{#100}}%

\newcommand{\nmono}{\nmonov{50}}%

\newcommand{\Nmonov}[2]{\Vcase{\NMONO}{#1}{#200}}%

\newcommand{\Nmono}[1]{\Nmonov{#1}{50}}%

\newcommand{\nmonOv}[2]{\vcasE{\NMONO}{#1}{#200}}%

\newcommand{\nmonO}[1]{\nmonOv{#1}{50}}%

\newcommand{\nepiv}[1]{\vcase{\NEPI}{#100}}%

\newcommand{\nepi}{\nepiv{50}}%

\newcommand{\Nepiv}[2]{\Vcase{\NEPI}{#1}{#200}}%

\newcommand{\Nepi}[1]{\Nepiv{#1}{50}}%

\newcommand{\nepIv}[2]{\vcasE{\NEPI}{#1}{#200}}%

\newcommand{\nepI}[1]{\nepIv{#1}{50}}%

\newcommand{\nbimov}[1]{\vcase{\NBIMO}{#100}}%

\newcommand{\nbimo}{\nbimov{50}}%

\newcommand{\Nbimov}[2]{\Vcase{\NBIMO}{#1}{#200}}%

\newcommand{\Nbimo}[1]{\Nbimov{#1}{50}}%

\newcommand{\nbimOv}[2]{\vcasE{\NBIMO}{#1}{#200}}%

\newcommand{\nbimO}[1]{\nbimOv{#1}{50}}%

\newcommand{\nisov}[1]{\vcasE{\NAR}{\cong}{#100}}%

\newcommand{\niso}{\nisov{50}}%

\newcommand{\Nisov}[2]%
{\Vbicase{\NAR}{#1\hspace{-2pt}}{\hspace{-2pt}\cong}{#200}}%

\newcommand{\Niso}[1]{\Nisov{#1}{50}}%

\newcommand{\nbiarv}[1]{\vbicase{\NBIAR}{#100}}%

\newcommand{\nbiar}{\nbiarv{50}}%

\newcommand{\Nbiarv}[3]{\Vbicase{\NBIAR}{#1}{#2}{#300}}%

\newcommand{\Nbiar}[2]{\Nbiarv{#1}{#2}{50}}%

\newcommand{\neqlv}[1]{\vbicase{\SEQL}{#100}}%

\newcommand{\neql}{\neqlv{50}}%

\newcommand{\nadjarv}[1]{\vbicase{\NADJAR}{#100}}%

\newcommand{\nadjar}{\nadjarv{50}}%

\newcommand{\Nadjarv}[3]{\Vbicase{\NADJAR}{#1}{#2}{#300}}%

\newcommand{\Nadjar}[2]{\Nadjarv{#1}{#2}{50}}%


\newcommand{\fdcase}[3]{\begin{picture}(0,0)%
\put(0,-150){#1}%
\truex{200}\truey{600}\truez{600}%
\put(-\value{x},-\value{x}){\makebox(0,\value{z})[r]{${#2}$}}%
\put(\value{x},-\value{y}){\makebox(0,\value{z})[l]{${#3}$}}%
\end{picture}}%

\newcommand{\fdbicase}[3]{\begin{picture}(0,0)%
\put(0,-150){#1}%
\truex{800}\truey{50}%
\put(-\value{x},\value{y}){${#2}$}%
\truex{200}\truey{950}%
\put(\value{x},-\value{y}){${#3}$}%
\end{picture}}%


\newcommand{\NEAR}{\begin{picture}(0,0)%
\put(-2900,-2900){\vector(1,1){5800}}%
\end{picture}}%

\newcommand{\NEDOTAR}%
{\truex{100}\truey{212}%
\NUMBEROFDOTS=5800%
\divide\NUMBEROFDOTS by \value{y}%
\begin{picture}(0,0)%
\multiput(-2900,-2900)(\value{y},\value{y}){\NUMBEROFDOTS}%
{\circle*{\value{x}}}%
\put(2900,2900){\vector(1,1){0}}%
\end{picture}}%

\newcommand{\NEMONO}{\begin{picture}(0,0)%
\put(-2900,-2900){\vector(1,1){5800}}%
\put(-2900,-2900){\begin{picture}(0,0)%
\truex{141}%
\put(\value{x},\value{x}){\makebox(0,0){$\times$}}%
\end{picture}}\end{picture}}%

\newcommand{\NEEPI}{\begin{picture}(0,0)%
\put(-2900,-2900){\vector(1,1){5800}}%
\put(2900,2900){\begin{picture}(0,0)%
\truex{545}%
\put(-\value{x},-\value{x}){\makebox(0,0){$\times$}}%
\end{picture}}\end{picture}}%

\newcommand{\NEBIMO}{\begin{picture}(0,0)%
\put(-2900,-2900){\vector(1,1){5800}}%
\put(2900,2900){\begin{picture}(0,0)%
\truex{545}%
\put(-\value{x},-\value{x}){\makebox(0,0){$\times$}}%
\end{picture}}
\put(-2900,-2900){\begin{picture}(0,0)%
\truex{141}%
\put(\value{x},\value{x}){\makebox(0,0){$\times$}}%
\end{picture}}\end{picture}}%

\newcommand{\NEBIAR}{\begin{picture}(0,0)%
\put(-2900,-2900){\begin{picture}(0,0)%
\truex{141}%
\put(-\value{x},\value{x}){\vector(1,1){5800}}%
\put(\value{x},-\value{x}){\vector(1,1){5800}}%
\end{picture}}\end{picture}}%

\newcommand{\NEEQL}{\begin{picture}(0,0)%
\put(-2900,-2900){\begin{picture}(0,0)%
\truex{70}%
\put(-\value{x},\value{x}){\line(1,1){5800}}%
\put(\value{x},-\value{x}){\line(1,1){5800}}%
\end{picture}}\end{picture}}%

\newcommand{\NEADJAR}{\begin{picture}(0,0)%
\put(-2900,-2900){\begin{picture}(0,0)%
\truex{141}%
\put(\value{x},-\value{x}){\vector(1,1){5800}}%
\end{picture}}%
\put(2900,2900){\begin{picture}(0,0)%
\truex{141}%
\put(-\value{x},\value{x}){\vector(-1,-1){5800}}%
\end{picture}}\end{picture}}%

\newcommand{\NEARV}[1]{\begin{picture}(0,0)%
\put(0,0){\makebox(0,0){\begin{picture}(#1,#1)%
\put(0,0){\vector(1,1){#1}}\end{picture}}}%
\end{picture}}%


\newcommand{\near}{\fdcase{\NEAR}{}{}}%

\newcommand{\Near}[1]{\fdcase{\NEAR}{#1}{}}%

\newcommand{\neaR}[1]{\fdcase{\NEAR}{}{#1}}%

\newcommand{\nedotar}{\fdcase{\NEDOTAR}{}{}}%

\newcommand{\Nedotar}[1]{\fdcase{\NEDOTAR}{#1}{}}%

\newcommand{\nedotaR}[1]{\fdcase{\NEDOTAR}{}{#1}}%

\newcommand{\nemono}{\fdcase{\NEMONO}{}{}}%

\newcommand{\Nemono}[1]{\fdcase{\NEMONO}{#1}{}}%

\newcommand{\nemonO}[1]{\fdcase{\NEMONO}{}{#1}}%

\newcommand{\neepi}{\fdcase{\NEEPI}{}{}}%

\newcommand{\Neepi}[1]{\fdcase{\NEEPI}{#1}{}}%

\newcommand{\neepI}[1]{\fdcase{\NEEPI}{}{#1}}%

\newcommand{\nebimo}{\fdcase{\NEBIMO}{}{}}%

\newcommand{\Nebimo}[1]{\fdcase{\NEBIMO}{#1}{}}%

\newcommand{\nebimO}[1]{\fdcase{\NEBIMO}{}{#1}}%

\newcommand{\neiso}{\fdcase{\NEAR}{\hspace{-2pt}\cong}{}}%

\newcommand{\Neiso}[1]{\fdcase{\NEAR}{\hspace{-2pt}\cong}{#1}}%

\newcommand{\nebiar}{\fdbicase{\NEBIAR}{}{}}%

\newcommand{\Nebiar}[2]{\fdbicase{\NEBIAR}{#1}{#2}}%

\newcommand{\neeql}{\fdbicase{\NEEQL}{}{}}%

\newcommand{\neadjar}{\fdbicase{\NEADJAR}{}{}}%

\newcommand{\Neadjar}[2]{\fdbicase{\NEADJAR}{#1}{#2}}%


\newcommand{\nearv}[1]{\fdcase{\NEARV{#100}}{}{}}%

\newcommand{\Nearv}[2]{\fdcase{\NEARV{#200}}{#1}{}}%

\newcommand{\neaRv}[2]{\fdcase{\NEARV{#200}}{}{#1}}%


\newcommand{\SWAR}{\begin{picture}(0,0)%
\put(2900,2900){\vector(-1,-1){5800}}%
\end{picture}}%

\newcommand{\SWDOTAR}%
{\truex{100}\truey{212}%
\NUMBEROFDOTS=5800%
\divide\NUMBEROFDOTS by \value{y}%
\begin{picture}(0,0)%
\multiput(2900,2900)(-\value{y},-\value{y}){\NUMBEROFDOTS}%
{\circle*{\value{x}}}%
\put(-2900,-2900){\vector(-1,-1){0}}%
\end{picture}}%

\newcommand{\SWMONO}{\begin{picture}(0,0)%
\put(2900,2900){\vector(-1,-1){5800}}%
\put(2900,2900){\begin{picture}(0,0)%
\truex{141}%
\put(-\value{x},-\value{x}){\makebox(0,0){$\times$}}%
\end{picture}}\end{picture}}%

\newcommand{\SWEPI}{\begin{picture}(0,0)%
\put(2900,2900){\vector(-1,-1){5800}}%
\put(-2900,-2900){\begin{picture}(0,0)%
\truex{525}%
\put(\value{x},\value{x}){\makebox(0,0){$\times$}}%
\end{picture}}\end{picture}}%

\newcommand{\SWBIMO}{\begin{picture}(0,0)%
\put(2900,2900){\vector(-1,-1){5800}}%
\put(2900,2900){\begin{picture}(0,0)%
\truex{141}%
\put(-\value{x},-\value{x}){\makebox(0,0){$\times$}}%
\end{picture}}%
\put(-2900,-2900){\begin{picture}(0,0)%
\truex{525}%
\put(\value{x},\value{x}){\makebox(0,0){$\times$}}%
\end{picture}}\end{picture}}%

\newcommand{\SWBIAR}{\begin{picture}(0,0)%
\put(2900,2900){\begin{picture}(0,0)%
\truex{141}%
\put(\value{x},-\value{x}){\vector(-1,-1){5800}}%
\put(-\value{x},\value{x}){\vector(-1,-1){5800}}%
\end{picture}}\end{picture}}%

\newcommand{\SWADJAR}{\begin{picture}(0,0)%
\put(-2900,-2900){\begin{picture}(0,0)%
\truex{141}%
\put(-\value{x},\value{x}){\vector(1,1){5800}}%
\end{picture}}%
\put(2900,2900){\begin{picture}(0,0)%
\truex{141}%
\put(\value{x},-\value{x}){\vector(-1,-1){5800}}%
\end{picture}}\end{picture}}%

\newcommand{\SWARV}[1]{\begin{picture}(0,0)%
\put(0,0){\makebox(0,0){\begin{picture}(#1,#1)%
\put(#1,#1){\vector(-1,-1){#1}}\end{picture}}}%
\end{picture}}%


\newcommand{\swar}{\fdcase{\SWAR}{}{}}%

\newcommand{\Swar}[1]{\fdcase{\SWAR}{#1}{}}%

\newcommand{\swaR}[1]{\fdcase{\SWAR}{}{#1}}%

\newcommand{\swdotar}{\fdcase{\SWDOTAR}{}{}}%

\newcommand{\Swdotar}[1]{\fdcase{\SWDOTAR}{#1}{}}%

\newcommand{\swdotaR}[1]{\fdcase{\SWDOTAR}{}{#1}}%

\newcommand{\swmono}{\fdcase{\SWMONO}{}{}}%

\newcommand{\Swmono}[1]{\fdcase{\SWMONO}{#1}{}}%

\newcommand{\swmonO}[1]{\fdcase{\SWMONO}{}{#1}}%

\newcommand{\swepi}{\fdcase{\SWEPI}{}{}}%

\newcommand{\Swepi}[1]{\fdcase{\SWEPI}{#1}{}}%

\newcommand{\swepI}[1]{\fdcase{\SWEPI}{}{#1}}%

\newcommand{\swbimo}{\fdcase{\SWBIMO}{}{}}%

\newcommand{\Swbimo}[1]{\fdcase{\SWBIMO}{#1}{}}%

\newcommand{\swbimO}[1]{\fdcase{\SWBIMO}{}{#1}}%

\newcommand{\swiso}{\fdcase{\SWAR}{\hspace{-2pt}\cong}{}}%

\newcommand{\Swiso}[1]{\fdcase{\SWAR}{\hspace{-2pt}\cong}{#1}}%

\newcommand{\swbiar}{\fdbicase{\SWBIAR}{}{}}%

\newcommand{\Swbiar}[2]{\fdbicase{\SWBIAR}{#1}{#2}}%

\newcommand{\sweql}{\fdbicase{\NEEQL}{}{}}%

\newcommand{\swadjar}{\fdbicase{\SWADJAR}{}{}}%

\newcommand{\Swadjar}[2]{\fdbicase{\SWADJAR}{#1}{#2}}%


\newcommand{\swarv}[1]{\fdcase{\SWARV{#100}}{}{}}%

\newcommand{\Swarv}[2]{\fdcase{\SWARV{#200}}{#1}{}}%

\newcommand{\swaRv}[2]{\fdcase{\SWARV{#200}}{}{#1}}%


\newcommand{\sdcase}[3]{\begin{picture}(0,0)%
\put(0,-150){#1}%
\truex{100}\truez{600}%
\put(\value{x},\value{x}){\makebox(0,\value{z})[l]{${#2}$}}%
\truex{300}\truey{800}%
\put(-\value{x},-\value{y}){\makebox(0,\value{z})[r]{${#3}$}}%
\end{picture}}%

\newcommand{\sdbicase}[3]{\begin{picture}(0,0)%
\put(0,-150){#1}%
\truex{250}\truey{600}\truez{850}%
\put(\value{x},\value{x}){\makebox(0,\value{y})[l]{${#2}$}}%
\put(-\value{x},-\value{z}){\makebox(0,\value{y})[r]{${#3}$}}%
\end{picture}}%


\newcommand{\SEAR}{\begin{picture}(0,0)%
\put(-2900,2900){\vector(1,-1){5800}}%
\end{picture}}%

\newcommand{\SEDOTAR}%
{\truex{100}\truey{212}%
\NUMBEROFDOTS=5800%
\divide\NUMBEROFDOTS by \value{y}%
\begin{picture}(0,0)%
\multiput(-2900,2900)(\value{y},-\value{y}){\NUMBEROFDOTS}%
{\circle*{\value{x}}}%
\put(2900,-2900){\vector(1,-1){0}}%
\end{picture}}%

\newcommand{\SEMONO}{\begin{picture}(0,0)%
\put(-2900,2900){\vector(1,-1){5800}}%
\put(-2900,2900){\begin{picture}(0,0)%
\truex{141}%
\put(\value{x},-\value{x}){\makebox(0,0){$\times$}}%
\end{picture}}\end{picture}}%

\newcommand{\SEEPI}{\begin{picture}(0,0)%
\put(-2900,2900){\vector(1,-1){5800}}%
\put(2900,-2900){\begin{picture}(0,0)%
\truex{525}%
\put(-\value{x},\value{x}){\makebox(0,0){$\times$}}%
\end{picture}}\end{picture}}%

\newcommand{\SEBIMO}{\begin{picture}(0,0)%
\put(-2900,2900){\vector(1,-1){5800}}%
\put(-2900,2900){\begin{picture}(0,0)%
\truex{141}%
\put(\value{x},-\value{x}){\makebox(0,0){$\times$}}%
\end{picture}}%
\put(2900,-2900){\begin{picture}(0,0)%
\truex{525}%
\put(-\value{x},\value{x}){\makebox(0,0){$\times$}}%
\end{picture}}\end{picture}}%

\newcommand{\SEBIAR}{\begin{picture}(0,0)%
\put(-2900,2900){\begin{picture}(0,0)%
\truex{141}
\put(-\value{x},-\value{x}){\vector(1,-1){5800}}%
\put(\value{x},\value{x}){\vector(1,-1){5800}}%
\end{picture}}\end{picture}}%

\newcommand{\SEEQL}{\begin{picture}(0,0)%
\put(-2900,2900){\begin{picture}(0,0)%
\truex{70}%
\put(-\value{x},-\value{x}){\line(1,-1){5800}}%
\put(\value{x},\value{x}){\line(1,-1){5800}}%
\end{picture}}\end{picture}}%

\newcommand{\SEADJAR}{\begin{picture}(0,0)%
\put(-2900,2900){\begin{picture}(0,0)%
\truex{141}%
\put(-\value{x},-\value{x}){\vector(1,-1){5800}}%
\end{picture}}%
\put(2900,-2900){\begin{picture}(0,0)%
\truex{141}%
\put(\value{x},\value{x}){\vector(-1,1){5800}}%
\end{picture}}\end{picture}}%

\newcommand{\SEARV}[1]{\begin{picture}(0,0)%
\put(0,0){\makebox(0,0){\begin{picture}(#1,#1)%
\put(0,#1){\vector(1,-1){#1}}\end{picture}}}%
\end{picture}}%


\newcommand{\sear}{\sdcase{\SEAR}{}{}}%

\newcommand{\Sear}[1]{\sdcase{\SEAR}{#1}{}}%

\newcommand{\seaR}[1]{\sdcase{\SEAR}{}{#1}}%

\newcommand{\sedotar}{\sdcase{\SEDOTAR}{}{}}%

\newcommand{\Sedotar}[1]{\sdcase{\SEDOTAR}{#1}{}}%

\newcommand{\sedotaR}[1]{\sdcase{\SEDOTAR}{}{#1}}%

\newcommand{\semono}{\sdcase{\SEMONO}{}{}}%

\newcommand{\Semono}[1]{\sdcase{\SEMONO}{#1}{}}%

\newcommand{\semonO}[1]{\sdcase{\SEMONO}{}{#1}}%

\newcommand{\seepi}{\sdcase{\SEEPI}{}{}}%

\newcommand{\Seepi}[1]{\sdcase{\SEEPI}{#1}{}}%

\newcommand{\seepI}[1]{\sdcase{\SEEPI}{}{#1}}%

\newcommand{\sebimo}{\sdcase{\SEBIMO}{}{}}%

\newcommand{\Sebimo}[1]{\sdcase{\SEBIMO}{#1}{}}%

\newcommand{\sebimO}[1]{\sdcase{\SEBIMO}{}{#1}}%

\newcommand{\seiso}{\sdcase{\SEAR}{\hspace{-2pt}\cong}{}}%

\newcommand{\Seiso}[1]{\sdcase{\SEAR}{\hspace{-2pt}\cong}{#1}}%

\newcommand{\sebiar}{\sdbicase{\SEBIAR}{}{}}%

\newcommand{\Sebiar}[2]{\sdbicase{\SEBIAR}{#1}{#2}}%

\newcommand{\seeql}{\sdbicase{\SEEQL}{}{}}%

\newcommand{\seadjar}{\sdbicase{\SEADJAR}{}{}}%

\newcommand{\Seadjar}[2]{\sdbicase{\SEADJAR}{#1}{#2}}%


\newcommand{\searv}[1]{\sdcase{\SEARV{#100}}{}{}}%

\newcommand{\Searv}[2]{\sdcase{\SEARV{#200}}{#1}{}}%

\newcommand{\seaRv}[2]{\sdcase{\SEARV{#200}}{}{#1}}%


\newcommand{\NWAR}{\begin{picture}(0,0)%
\put(2900,-2900){\vector(-1,1){5800}}%
\end{picture}}%

\newcommand{\NWDOTAR}%
{\truex{100}\truey{212}%
\NUMBEROFDOTS=5800%
\divide\NUMBEROFDOTS by \value{y}%
\begin{picture}(0,0)%
\multiput(2900,-2900)(-\value{y},\value{y}){\NUMBEROFDOTS}%
{\circle*{\value{x}}}%
\put(-2900,2900){\vector(-1,1){0}}%
\end{picture}}%

\newcommand{\NWMONO}{\begin{picture}(0,0)%
\put(2900,-2900){\vector(-1,1){5800}}%
\put(2900,-2900){\begin{picture}(0,0)%
\truex{141}%
\put(-\value{x},\value{x}){\makebox(0,0){$\times$}}%
\end{picture}}\end{picture}}%

\newcommand{\NWEPI}{\begin{picture}(0,0)%
\put(2900,-2900){\vector(-1,1){5800}}%
\put(-2900,2900){\begin{picture}(0,0)%
\truex{525}%
\put(\value{x},-\value{x}){\makebox(0,0){$\times$}}%
\end{picture}}\end{picture}}%

\newcommand{\NWBIMO}{\begin{picture}(0,0)%
\put(2900,-2900){\vector(-1,1){5800}}%
\put(2900,-2900){\begin{picture}(0,0)%
\truex{141}%
\put(-\value{x},\value{x}){\makebox(0,0){$\times$}}%
\end{picture}}%
\put(-2900,2900){\begin{picture}(0,0)%
\truex{525}%
\put(\value{x},-\value{x}){\makebox(0,0){$\times$}}%
\end{picture}}\end{picture}}%

\newcommand{\NWBIAR}{\begin{picture}(0,0)%
\put(2900,-2900){\begin{picture}(0,0)%
\truex{141}%
\put(\value{x},\value{x}){\vector(-1,1){5800}}%
\end{picture}}%
\put(2900,-2900){\begin{picture}(0,0)%
\truex{141}
\put(-\value{x},-\value{x}){\vector(-1,1){5800}}%
\end{picture}}\end{picture}}%

\newcommand{\NWADJAR}{\begin{picture}(0,0)%
\put(-2900,2900){\begin{picture}(0,0)%
\truex{141}%
\put(\value{x},\value{x}){\vector(1,-1){5800}}%
\end{picture}}%
\put(2900,-2900){\begin{picture}(0,0)%
\truex{141}%
\put(-\value{x},-\value{x}){\vector(-1,1){5800}}%
\end{picture}}\end{picture}}%

\newcommand{\NWARV}[1]{\begin{picture}(0,0)%
\put(0,0){\makebox(0,0){\begin{picture}(#1,#1)%
\put(#1,0){\vector(-1,1){#1}}\end{picture}}}%
\end{picture}}%


\newcommand{\nwar}{\sdcase{\NWAR}{}{}}%

\newcommand{\Nwar}[1]{\sdcase{\NWAR}{#1}{}}%

\newcommand{\nwaR}[1]{\sdcase{\NWAR}{}{#1}}%

\newcommand{\nwdotar}{\sdcase{\NWDOTAR}{}{}}%

\newcommand{\Nwdotar}[1]{\sdcase{\NWDOTAR}{#1}{}}%

\newcommand{\nwdotaR}[1]{\sdcase{\NWDOTAR}{}{#1}}%

\newcommand{\nwmono}{\sdcase{\NWMONO}{}{}}%

\newcommand{\Nwmono}[1]{\sdcase{\NWMONO}{#1}{}}%

\newcommand{\nwmonO}[1]{\sdcase{\NWMONO}{}{#1}}%

\newcommand{\nwepi}{\sdcase{\NWEPI}{}{}}%

\newcommand{\Nwepi}[1]{\sdcase{\NWEPI}{#1}{}}%

\newcommand{\nwepI}[1]{\sdcase{\NWEPI}{}{#1}}%

\newcommand{\nwbimo}{\sdcase{\NWBIMO}{}{}}%

\newcommand{\Nwbimo}[1]{\sdcase{\NWBIMO}{#1}{}}%

\newcommand{\nwbimO}[1]{\sdcase{\NWBIMO}{}{#1}}%

\newcommand{\nwiso}{\sdcase{\NWAR}{\hspace{-2pt}\cong}{}}%

\newcommand{\Nwiso}[1]{\sdcase{\NWAR}{\hspace{-2pt}\cong}{#1}}%

\newcommand{\nwbiar}{\sdbicase{\NWBIAR}{}{}}%

\newcommand{\Nwbiar}[2]{\sdbicase{\NWBIAR}{#1}{#2}}%

\newcommand{\nweql}{\sdbicase{\SEEQL}{}{}}%

\newcommand{\nwadjar}{\sdbicase{\NWADJAR}{}{}}%

\newcommand{\Nwadjar}[2]{\sdbicase{\NWADJAR}{#1}{#2}}%


\newcommand{\nwarv}[1]{\sdcase{\NWARV{#100}}{}{}}%

\newcommand{\Nwarv}[2]{\sdcase{\NWARV{#200}}{#1}{}}%

\newcommand{\nwaRv}[2]{\sdcase{\NWARV{#200}}{}{#1}}%



\newcommand{\ENEAR}[2]%
{\makebox[0pt]{\begin{picture}(0,0)%
\put(0,-150){\makebox(0,0){\begin{picture}(0,0)%
\put(-6600,-3300){\vector(2,1){13200}}%
\truex{200}\truey{800}\truez{600}%
\put(-\value{x},\value{x}){\makebox(0,\value{z})[r]{${#1}$}}%
\put(\value{x},-\value{y}){\makebox(0,\value{z})[l]{${#2}$}}%
\end{picture}}}\end{picture}}}%

\newcommand{\enear}{\ENEAR{}{}}%

\newcommand{\Enear}[1]{\ENEAR{#1}{}}%

\newcommand{\eneaR}[1]{\ENEAR{}{#1}}%

\newcommand{\ESEAR}[2]%
{\makebox[0pt]{\begin{picture}(0,0)%
\put(0,-150){\makebox(0,0){\begin{picture}(0,0)%
\put(-6600,3300){\vector(2,-1){13200}}%
\truex{200}\truey{800}\truez{600}%
\put(\value{x},\value{x}){\makebox(0,\value{z})[l]{${#1}$}}%
\put(-\value{x},-\value{y}){\makebox(0,\value{z})[r]{${#2}$}}%
\end{picture}}}\end{picture}}}%

\newcommand{\esear}{\ESEAR{}{}}%

\newcommand{\Esear}[1]{\ESEAR{#1}{}}%

\newcommand{\eseaR}[1]{\ESEAR{}{#1}}%

\newcommand{\WNWAR}[2]%
{\makebox[0pt]{\begin{picture}(0,0)%
\put(0,-150){\makebox(0,0){\begin{picture}(0,0)%
\put(6600,-3300){\vector(-2,1){13200}}%
\truex{200}\truey{800}\truez{600}%
\put(\value{x},\value{x}){\makebox(0,\value{z})[l]{${#1}$}}%
\put(-\value{x},-\value{y}){\makebox(0,\value{z})[r]{${#2}$}}%
\end{picture}}}\end{picture}}}%

\newcommand{\wnwar}{\WNWAR{}{}}%

\newcommand{\Wnwar}[1]{\WNWAR{#1}{}}%

\newcommand{\wnwaR}[1]{\WNWAR{}{#1}}%

\newcommand{\WSWAR}[2]%
{\makebox[0pt]{\begin{picture}(0,0)%
\put(0,-150){\makebox(0,0){\begin{picture}(0,0)%
\put(6600,3300){\vector(-2,-1){13200}}%
\truex{200}\truey{800}\truez{600}%
\put(-\value{x},\value{x}){\makebox(0,\value{z})[r]{${#1}$}}%
\put(\value{x},-\value{y}){\makebox(0,\value{z})[l]{${#2}$}}%
\end{picture}}}\end{picture}}}%

\newcommand{\wswar}{\WSWAR{}{}}%

\newcommand{\Wswar}[1]{\WSWAR{#1}{}}%

\newcommand{\wswaR}[1]{\WSWAR{}{#1}}%



\newcommand{\NNEAR}[2]%
{\raisebox{-1pt}[0pt][0pt]{\begin{picture}(0,0)%
\put(0,0){\makebox(0,0){\begin{picture}(0,0)%
\put(-3300,-6600){\vector(1,2){6600}}%
\truex{100}\truez{600}%
\put(-\value{x},\value{x}){\makebox(0,\value{z})[r]{${#1}$}}%
\put(\value{x},-\value{z}){\makebox(0,\value{z})[l]{${#2}$}}%
\end{picture}}}\end{picture}}}%

\newcommand{\nnear}{\NNEAR{}{}}%

\newcommand{\Nnear}[1]{\NNEAR{#1}{}}%

\newcommand{\nneaR}[1]{\NNEAR{}{#1}}%

\newcommand{\SSWAR}[2]%
{\raisebox{-1pt}[0pt][0pt]{\begin{picture}(0,0)%
\put(0,0){\makebox(0,0){\begin{picture}(0,0)%
\put(3300,6600){\vector(-1,-2){6600}}%
\truex{100}\truez{600}%
\put(-\value{x},\value{x}){\makebox(0,\value{z})[r]{${#1}$}}%
\put(\value{x},-\value{z}){\makebox(0,\value{z})[l]{${#2}$}}%
\end{picture}}}\end{picture}}}%

\newcommand{\sswar}{\SSWAR{}{}}%

\newcommand{\Sswar}[1]{\SSWAR{#1}{}}%

\newcommand{\sswaR}[1]{\SSWAR{}{#1}}%

\newcommand{\SSEAR}[2]%
{\raisebox{-1pt}[0pt][0pt]{\begin{picture}(0,0)%
\put(0,0){\makebox(0,0){\begin{picture}(0,0)%
\put(-3300,6600){\vector(1,-2){6600}}%
\truex{200}\truez{600}%
\put(\value{x},\value{x}){\makebox(0,\value{z})[l]{${#1}$}}%
\put(-\value{x},-\value{z}){\makebox(0,\value{z})[r]{${#2}$}}%
\end{picture}}}\end{picture}}}%

\newcommand{\ssear}{\SSEAR{}{}}%

\newcommand{\Ssear}[1]{\SSEAR{#1}{}}%

\newcommand{\sseaR}[1]{\SSEAR{}{#1}}%

\newcommand{\NNWAR}[2]%
{\raisebox{-1pt}[0pt][0pt]{\begin{picture}(0,0)%
\put(0,0){\makebox(0,0){\begin{picture}(0,0)%
\put(3300,-6600){\vector(-1,2){6600}}%
\truex{200}\truez{600}%
\put(\value{x},\value{x}){\makebox(0,\value{z})[l]{${#1}$}}%
\put(-\value{x},-\value{z}){\makebox(0,\value{z})[r]{${#2}$}}%
\end{picture}}}\end{picture}}}%

\newcommand{\nnwar}{\NNWAR{}{}}%

\newcommand{\Nnwar}[1]{\NNWAR{#1}{}}%

\newcommand{\nnwaR}[1]{\NNWAR{}{#1}}%



\newcommand{\Necurve}[2]%
{\begin{picture}(0,0)%
\truex{1300}\truey{2000}\truez{200}%
\put(0,\value{x}){\oval(#200,\value{y})[t]}%
\put(0,\value{x}){\makebox(0,0){\begin{picture}(#200,0)%
\put(#200,0){\vector(0,-1){\value{z}}}%
\put(0,0){\line(0,-1){\value{z}}}\end{picture}}}%
\truex{2500}%
\put(0,\value{x}){\makebox(0,0)[b]{${#1}$}}%
\end{picture}}%

\newcommand{\necurve}[1]{\Necurve{}{#1}}%

\newcommand{\Nwcurve}[2]%
{\begin{picture}(0,0)%
\truex{1300}\truey{2000}\truez{200}%
\put(0,\value{x}){\oval(#200,\value{y})[t]}%
\put(0,\value{x}){\makebox(0,0){\begin{picture}(#200,0)%
\put(#200,0){\line(0,-1){\value{z}}}%
\put(0,0){\vector(0,-1){\value{z}}}\end{picture}}}%
\truex{2500}%
\put(0,\value{x}){\makebox(0,0)[b]{${#1}$}}%
\end{picture}}%

\newcommand{\nwcurve}[1]{\Nwcurve{}{#1}}%

\newcommand{\Securve}[2]%
{\begin{picture}(0,0)%
\truex{1300}\truey{2000}\truez{200}%
\put(0,-\value{x}){\oval(#200,\value{y})[b]}%
\put(0,-\value{x}){\makebox(0,0){\begin{picture}(#200,0)%
\put(#200,0){\vector(0,1){\value{z}}}%
\put(0,0){\line(0,1){\value{z}}}\end{picture}}}%
\truex{2500}%
\put(0,-\value{x}){\makebox(0,0)[t]{${#1}$}}%
\end{picture}}%

\newcommand{\securve}[1]{\Securve{}{#1}}%

\newcommand{\Swcurve}[2]%
{\begin{picture}(0,0)%
\truex{1300}\truey{2000}\truez{200}%
\put(0,-\value{x}){\oval(#200,\value{y})[b]}%
\put(0,-\value{x}){\makebox(0,0){\begin{picture}(#200,0)%
\put(#200,0){\line(0,1){\value{z}}}%
\put(0,0){\vector(0,1){\value{z}}}\end{picture}}}%
\truex{2500}%
\put(0,-\value{x}){\makebox(0,0)[t]{${#1}$}}%
\end{picture}}%

\newcommand{\swcurve}[1]{\Swcurve{}{#1}}%



\newcommand{\Escurve}[2]%
{\begin{picture}(0,0)%
\truex{1400}\truey{2000}\truez{200}%
\put(\value{x},0){\oval(\value{y},#200)[r]}%
\put(\value{x},0){\makebox(0,0){\begin{picture}(0,#200)%
\put(0,0){\vector(-1,0){\value{z}}}%
\put(0,#200){\line(-1,0){\value{z}}}\end{picture}}}%
\truex{2500}%
\put(\value{x},0){\makebox(0,0)[l]{${#1}$}}%
\end{picture}}%

\newcommand{\escurve}[1]{\Escurve{}{#1}}%

\newcommand{\Encurve}[2]%
{\begin{picture}(0,0)%
\truex{1400}\truey{2000}\truez{200}%
\put(\value{x},0){\oval(\value{y},#200)[r]}%
\put(\value{x},0){\makebox(0,0){\begin{picture}(0,#200)%
\put(0,0){\line(-1,0){\value{z}}}%
\put(0,#200){\vector(-1,0){\value{z}}}\end{picture}}}%
\truex{2500}%
\put(\value{x},0){\makebox(0,0)[l]{${#1}$}}%
\end{picture}}%

\newcommand{\encurve}[1]{\Encurve{}{#1}}%

\newcommand{\Wscurve}[2]%
{\begin{picture}(0,0)%
\truex{1300}\truey{2000}\truez{200}%
\put(-\value{x},0){\oval(\value{y},#200)[l]}%
\put(-\value{x},0){\makebox(0,0){\begin{picture}(0,#200)%
\put(0,0){\vector(1,0){\value{z}}}%
\put(0,#200){\line(1,0){\value{z}}}\end{picture}}}%
\truex{2400}%
\put(-\value{x},0){\makebox(0,0)[r]{${#1}$}}%
\end{picture}}%

\newcommand{\wscurve}[1]{\Wscurve{}{#1}}%

\newcommand{\Wncurve}[2]%
{\begin{picture}(0,0)%
\truex{1300}\truey{2000}\truez{200}%
\put(-\value{x},0){\oval(\value{y},#200)[l]}%
\put(-\value{x},0){\makebox(0,0){\begin{picture}(0,#200)%
\put(0,0){\line(1,0){\value{z}}}%
\put(0,#200){\vector(1,0){\value{z}}}\end{picture}}}%
\truex{2400}%
\put(-\value{x},0){\makebox(0,0)[r]{${#1}$}}%
\end{picture}}%

\newcommand{\wncurve}[1]{\Wncurve{}{#1}}%



\newcommand{\Freear}[8]{\begin{picture}(0,0)%
\put(#400,#500){\vector(#6,#7){#800}}%
\truex{#200}\truey{#300}%
\put(\value{x},\value{y}){$#1$}\end{picture}}%

\newcommand{\freear}[5]{\Freear{}{0}{0}{#1}{#2}{#3}{#4}{#5}}%


\newcount\SCALE%

\newcount\NUMBER%

\newcount\LINE%

\newcount\COLUMN%

\newcount\WIDTH%

\newcount\SOURCE%

\newcount\ARROW%

\newcount\TARGET%

\newcount\ARROWLENGTH%

\newcount\NUMBEROFDOTS%

\newcounter{x}%

\newcounter{y}%

\newcounter{z}%

\newcounter{horizontal}%

\newcounter{vertical}%

\newskip\itemlength%

\newskip\firstitem%

\newskip\seconditem%

\newcommand{\printarrow}{}%


\newcommand{\truex}[1]{%
\NUMBER=#1%
\multiply\NUMBER by 100%
\divide\NUMBER by \SCALE%
\setcounter{x}{\NUMBER}}%

\newcommand{\truey}[1]{%
\NUMBER=#1%
\multiply\NUMBER by 100%
\divide\NUMBER by \SCALE%
\setcounter{y}{\NUMBER}}%

\newcommand{\truez}[1]{%
\NUMBER=#1%
\multiply\NUMBER by 100%
\divide\NUMBER by \SCALE%
\setcounter{z}{\NUMBER}}%

\newcommand{\changecounters}[1]{%
\SOURCE=\ARROW%
\ARROW=\TARGET%
\settowidth{\itemlength}{#1}%
\ifdim \itemlength > 2800\unitlength%
\addtolength{\itemlength}{-2800\unitlength}%
\TARGET=\itemlength%
\divide\TARGET by 1310%
\multiply\TARGET by 100%
\divide\TARGET by \SCALE%
\else%
\TARGET=0%
\fi%
\ARROWLENGTH=5000%
\advance\ARROWLENGTH by -\SOURCE%
\advance\ARROWLENGTH by -\TARGET%
\advance\SOURCE by -\TARGET}%

\newcommand{\initialize}[1]{%
\LINE=0%
\COLUMN=0%
\WIDTH=0%
\ARROW=0%
\TARGET=0%
\changecounters{#1}%
\renewcommand{\printarrow}{#1}%
\begin{center}%
\vspace{10pt}%
\begin{picture}(0,0)}%

\newcommand{\DIAG}[1]{%
\SCALE=100%
\setlength{\unitlength}{655sp}%
\initialize{\mbox{$#1$}}}%

\newcommand{\DIAGV}[2]{%
\SCALE=#1%
\setlength{\unitlength}{655sp}%
\multiply\unitlength by \SCALE%
\divide\unitlength by 100%
\initialize{\mbox{$#2$}}}%

\newcommand{\n}[1]{%
\changecounters{\mbox{$#1$}}%
\put(\COLUMN,\LINE){\makebox(0,0){\printarrow}}%
\thinlines%
\renewcommand{\printarrow}{\mbox{$#1$}}%
\advance\COLUMN by 4000}%

\newcommand{\nn}[1]{%
\put(\COLUMN,\LINE){\makebox(0,0){\printarrow}}%
\thinlines%
\ifnum \WIDTH < \COLUMN%
\WIDTH=\COLUMN%
\else%
\fi%
\advance\LINE by -4000%
\COLUMN=0%
\ARROW=0%
\TARGET=0%
\changecounters{\mbox{$#1$}}%
\renewcommand{\printarrow}{\mbox{$#1$}}}%

\newcommand{\conclude}{%
\put(\COLUMN,\LINE){\makebox(0,0){\printarrow}}%
\thinlines%
\ifnum \WIDTH < \COLUMN%
\WIDTH=\COLUMN%
\else%
\fi%
\setcounter{horizontal}{\WIDTH}%
\setcounter{vertical}{-\LINE}%
\end{picture}}%

\newcommand{\diag}{%
\conclude%
\raisebox{0pt}[0pt][\value{vertical}\unitlength]{}%
\hspace*{\value{horizontal}\unitlength}%
\vspace{10pt}%
\end{center}%
\setlength{\unitlength}{1pt}}%

\newcommand{\diagv}[3]{%
\conclude%
\NUMBER=#1%
\rule{0pt}{\NUMBER pt}%
\hspace*{-#2pt}%
\raisebox{0pt}[0pt][\value{vertical}\unitlength]{}%
\hspace*{\value{horizontal}\unitlength}
\NUMBER=#3%
\advance\NUMBER by 10%
\vspace*{\NUMBER pt}%
\end{center}%
\setlength{\unitlength}{1pt}}%

\newcommand{\N}[1]%
{\raisebox{0pt}[7pt][0pt]{$#1$}}%

\newcommand{\movename}[3]{%
\hspace{#2pt}%
\raisebox{#3pt}[5pt][2pt]{\raisebox{#3pt}{$#1$}}%
\hspace{-#2pt}}%

\newcommand{\movearrow}[3]{%
\makebox[0pt]{%
\hspace{#2pt}\hspace{#2pt}%
\raisebox{#3pt}[0pt][0pt]{\raisebox{#3pt}{$#1$}}}}%

\newcommand{\movevertex}[3]{%
\mbox{\hspace{#2pt}%
\raisebox{#3pt}{\raisebox{#3pt}{$#1$}}%
\hspace{-#2pt}}}%

\newcommand{\crosslength}[2]{%
\settowidth{\firstitem}{#1}%
\settowidth{\seconditem}{#2}%
\ifdim\firstitem < \seconditem%
\itemlength=\seconditem%
\else%
\itemlength=\firstitem%
\fi%
\divide\itemlength by 2%
\hspace{\itemlength}}%

\newcommand{\cross}[2]{%
\crosslength{\mbox{$#1$}}{\mbox{$#2$}}%
\begin{picture}(0,0)%
\put(0,0){\makebox(0,0){$#1$}}%
\thinlines%
\put(0,0){\makebox(0,0){$#2$}}%
\thinlines%
\end{picture}%
\crosslength{\mbox{$#1$}}{\mbox{$#2$}}}%

\newcommand{\B}{\thicklines}


\newcommand{\adj}{\begin{picture}(9,6)%
\put(1,3){\line(1,0){6}}\put(7,0){\line(0,1){6}}%
\end{picture}}%

\newcommand{\com}{\begin{picture}(12,8)%
\put(6,4){\oval(8,8)[b]}\put(6,4){\oval(8,8)[r]}%
\put(6,8){\vector(-1,0){2}}\end{picture}}%

\newcommand{\Nat}[3]{\raisebox{-2pt}%
{\begin{picture}(34,15)%
\put(2,10){\vector(1,0){30}}%
\put(2,0){\vector(1,0){30}}%
\put(13,2){$\Downarrow$}%
\put(20,3){$\scriptstyle{#2}$}%
\put(4,11){$\scriptstyle{#1}$}%
\put(4,1){$\scriptstyle{#3}$}%
\end{picture}}}%

\newcommand{\nat}{\raisebox{-2pt}%
{\begin{picture}(34,10)%
\put(2,10){\vector(1,0){30}}%
\put(2,0){\vector(1,0){30}}%
\put(13,2){$\Downarrow$}%
\end{picture}}}%

\newcommand{\Binat}[5]{\raisebox{-7.5pt}%
{\begin{picture}(34,25)%
\put(2,20){\vector(1,0){30}}%
\put(2,10){\vector(1,0){30}}%
\put(2,0){\vector(1,0){30}}%
\put(13,12){$\Downarrow$}%
\put(13,2){$\Downarrow$}%
\put(20,13){$\scriptstyle{#2}$}%
\put(20,3){$\scriptstyle{#4}$}%
\put(4,21){$\scriptstyle{#1}$}%
\put(4,11){$\scriptstyle{#3}$}%
\put(4,1){$\scriptstyle{#5}$}%
\end{picture}}}%

\newcommand{\binat}{\raisebox{-7.5pt}%
{\begin{picture}(34,20)%
\put(2,20){\vector(1,0){30}}%
\put(2,10){\vector(1,0){30}}%
\put(2,0){\vector(1,0){30}}%
\put(13,12){$\Downarrow$}%
\put(13,2){$\Downarrow$}%
\end{picture}}}%

\setcounter{section}{-1}

 \title{Non K\"ahlerian surfaces with a cycle of rational curves}
\author {Georges Dloussky}
\date{06/2020}

\begin{abstract} Let $S$ be a compact complex surface in class VII$_0^+$ containing a cycle of rational curves $C=\sum D_j$. Let $D=C+A$ be the maximal connected divisor containing $C$. If there is another connected component of curves $C'$ then $C'$ is a cycle of rational curves, $A=0$ and $S$ is a Inoue-Hirzebruch surface. If there is only one connected component $D$ then  each connected component $A_i$ of $A$ is a chain of rational curves which intersects a curve $C_j$ of the cycle and for each curve $C_j$ of the cycle there at most one chain which meets $C_j$. In other words, we do not prove the existence of curves other those of the cycle $C$, but if some other curves exist the maximal divisor looks like the maximal divisor of a Kato surface with perhaps missing curves. The proof of this topological result is an application of Donaldson theorem on trivialization of the intersection form and of deformation theory. We apply this result to show that a twisted logarithmic $1$-form has a trivial vanishing divisor.
\end{abstract}

\maketitle

\tableofcontents

\section{Introduction}
A  minimal compact complex surface $S$ is said to be  of the class VII$_0$  of Kodaira if the first Betti number satisfies $b_1(S)=1$ and Kodaira dimension satisfies $\k(S)<0$. A surface $S$ is of class VII$_0^+$ if moreover 
$n:=b_2(S)>0$; these surfaces admit no nonconstant meromorphic functions. The cupproduct form $Q$ defines on $H^2(S,\bb R)$ a negative definite form and in 1987 S.K. Donaldson proved
\begin{Th}[\cite{DON}] If $X$ is  closed, oriented smooth $4$-manifold whose intersection form
$$Q:H^2(S,\bb Z)/{\rm Torsion}\to \bb Z$$
is negative definite, then the form is equivalent over the integers to the standard form $(-1)\oplus(-1)\oplus\cdots\oplus(-1)$.
\end{Th}
In other words, if $n=b_2(S)$, there are classes $e_0,\ldots,e_{n-1}$ such that $e_i.e_j=-\d_{ij}$. Under the further assumption that $S$ contains a cycle of rational curves, I. Nakamura \cite[Thm 1.5]{N2} proved the existence of these classes by deformation theory. In fact, $S$ can be deformed into a blown-up Hopf surface $n$ times \cite{N2} and the exceptional curves (or divisors) give the requested classes.
All known surfaces in class VII$_0^+$ are Kato surfaces i.e. contain Global Spherical Shells (GSS). The major problem in the classification of non-K\"ahlerian surfaces is to achieve the classification of surfaces $S$ of class VII$_0^+$. Are all surfaces in class VII$_0^+$ Kato surfaces ? Till now, the answer is positive for $b_2(S)=1$ by A. Teleman \cite{Te05}. By definition a GSS is a biholomorphic map $\f:U\to V$ from a neighbourhood $U\subset \bb C^2\setminus\{0\}$ of the sphere $S^3=\partial B^2$  onto an open set $V$ such that $\S=\f(S^3)$ does not disconnect $S$. By construction, rational curves are obtained as strict transform of exceptional curves. If for a subset $I\subset \{0,\ldots,n-1\}$ we denote $e_I:=\sum_{i\in I}e_i$, we can expect that the class of a curve is of the type $e_i-e_I$, $i\not\in I$. In \cite{D06} the following theorem is proved under, at one step, the existence of a Numerically AntiCanonical (NAC) divisor. In this article we avoid this assumption. Odd Inoue-Hirzebruch surfaces contain exactly one cycle are also called half Inoue-Hirzebruch surface. See \cite{D2} for a construction by contracting germs.
\begin{Th} Let $S$ be a minimal surface in class VII$_0^+$ with a cycle 
$$C=D_0+\cdots+D_{s-1} = -(e_r+\cdots+e_{n-1})$$
 of $s\ge 1$ rational curves. \\
1) $S$ is an odd Inoue-Hirzebruch surface if and only if $[H_1(C,\bb Z):H_1(S,\bb Z)]=2$. In this case $s=n$, $r=0$ and 
$$\sharp(C)-C^2=2b_2(S).$$

2) If $H_1(C,\bb Z)=H_1(S,\bb Z)$, then $s=r$,
$$\sharp(C)-C^2=b_2(S),$$
all curves are of type ${\bf a}$ i.e. $D_i = e_i-e_{I_i}$, $i\not\in I_i$, and $(I_i)_{0\le i\le s-1}$ is a partition of $[0,n-1]$.
\end{Th}
By this theorem it is possible to detect which exceptional curve of the first kind appears in a deformation of $S$ where one singular point of the cycle is smoothed (see theorem \ref{ejectioncourbeexc}). Since a surface with a cycle can be deformed into a blown-up primary Hopf surface \cite{N2} the groups $H_2(S,\bb Z)$ and $H^2(S,\bb Z)$ have no torsion.\\

By simple computations on the Donaldson classes it is easy to prove (see thm \ref{compconnexe}):
\begin{Th}\label{03}  Let $C=D_0+\cdots+D_{\a_i}+\cdots+D_{\a_{s-1}}$ be a cycle of $1\le s<n$ rational curves, $D=C+A$ be the maximal connected reduced divisor containing $C$.   Then there is a  subset $\mathfrak D\subset [0,n-1]$, such that 
$$D=C+A=-e_{\mathfrak D}\in H_2(S,\bb Z).$$
Moreover if $A$ is not empty, each connected component $A_i$ of $A$ called a tree is a chain of rational curves. Moreover two trees $A_j$, $A_k$ meet different curves of the cycle.
\end{Th}
In a series of papers \cite{Te05}, \cite{Te09}, \cite{Te10}, A. Teleman developped a strategy using gauge theory, successfull for $b_2(S)\le 3$, to show that any surface $S$ in class VII$_0^+$ contains a cycle of rational curves. The aim is to prove that a class of type $-e_I$ is represented by a one dimensional analytic subspace $C$. It is easy to show that with this condition $C$ has to contain a cycle. Theorem \ref{03} gives a reason (among others) why it is difficult to show the existence of the missing curves even when $b_2(S)=2$. Recall that if there are $b_2(S)$ rational curves the surface is Kato \cite{DOT3}.
Finally we obtain, with deformation theory, the following result which show with theorem \thesection.\ref{03} that the maximal divisor looks like the one of a Kato surface (see \cite{D1} or \cite{N1}).
\begin{Th} Let $S$ be a surface in class VII$_0^+$ containing a cycle $C$ of $s\ge 1$ rational curves and $D=C+A$ the connected component of the cycle. If there is another connected component $C'$ of curves, $C'$ is a cycle of rational curves, $A=0$ and $S$ is a Inoue-Hirzebruch surface.
\end{Th}
We apply this result to show that a twisted logarithmic $1$-form on a surface in class VII$_0^+$ cannot vanish along a divisor.
\section{Topological classes of curves and chains}
One suppose that $S$ contains a cycle of $s$ rational curves, $1\le s\le n$. When $s=1$ it is a rational curve with a double point, when $s=2$ it is the union of two rational curves $D_0,D_1$ such that $D_0.D_1=2$, and when $s\ge 3$, $D_0.D_1=D_1.D_2=\cdots=D_{s-1}.D_0=1$.\\
Since a surface in class VII$_0^+$ with a cycle of rational curves can be deformed into a blown-up Hopf surface, the groups  $H^2(S,\bb Z)$ and $H_2(S,\bb Z)$ have no torsion. By Poincaré duality they are isomorphic, therefore we shall denote in the same way a class $e_i$ in $H^2(S,\bb Z)$ and its dual in $H_2(S,\bb Z)$.

We recall basic facts used in \cite{N2}, \cite{D06} and in the next sections. For the sake of completeness we give a proof:

\begin{Lem} [\cite{N2} (2.5), (2.6), (2.7)] \label{N2(2.5)}
 1) Let $D$ be a nonsingular rational curve.  If  $D= a_0e_0+\cdots+a_{n-1}e_{n-1}\in H_2(S,\bb Z)$, then there exists a unique $i\in [0,n-1]$ such that $a_i=1\;{\rm or }\; -2$, and $a_j=0\;{\rm or}\; -1$ for $j\neq i$.\\
 2) Let $D_1$ and $D_2$ two distinct divisors without common irreducible component such that $D_1= e_{i_1}-e_{I_1}$, $i_1\not\in I_1$ and $D_2= e_{i_2}-e_{I_2}$, $i_2\not\in I_2$, then $i_1\neq i_2$.
 \end{Lem}
Proof: 1) For a rational curve we have $0=KD+D^2+2=2-\sum_k (a_k^2+a_k)$, therefore $a_i^2+a_i=2$ for a unique index and $a_k^2+a_k=0$ for the  others.\\
2) If $i_1=i_2$, $(e_{i_1}-e_{I_1}).(e_{i_2}-e_{I_2})=-1+e_{I_1}.e_{I_2}<0$ which is impossible.\hfill$\Box$
\begin{Lem}[\cite{N2} (2.4)] \label{N2(2.4)}
Let $S$ be a surface of class VII$_0^+$ with a cycle $C=D_0+\cdots+D_{s-1}= -(e_r+\cdots+e_{n-1})$ of rational curves and without divisor $D$ such that $D^2=0$. Let $L_{I}:=\bigotimes_{i\in I}L_{i}$, $L=L_{I}\otimes F$, $F\in H^1(S,\bb C^\star)$ for any subset $I\subset[0,n-1]$. Then we have:\\
i) If $I\neq [0,n-1]$, then $H^q(S,L)=0$ for any $q$.\\
ii) If $L\otimes \mathcal O_{C}=\mathcal O_{C}$, then $I=[0,r-1]$, 
and $F=\mathcal O_{S}$, $K_{S}\otimes L^\star\otimes[C]=\mathcal O_{S}$.\\
iii) If $LD_{i}=0$ for any irreducible component $D_{i}$ of $C$, then 
$I=[0,r-1]$.
\end{Lem}

Let $S$ be a surface containing a cycle $C=D_0+\cdots+D_{s-1}$ of $s$ rational curves. Lemma \thesection.\ref{N2(2.5)} shows that in the cycle there are two possible types of  nonsingular rational curves:
 \begin{itemize}
 \item (\bf Type a) $D_i= e_i-e_{I_i}$ with $i\not\in I_i$,
 \item (\bf Type b) $D_i=-2e_i-e_{I_i}$ with $i\not\in I_i$.
 \end{itemize}
 we shall show that the same dichotomy applies for chains of nonsingular rational curves.
\begin{Lem} \label{chaineb} Let $D'=D_j+\cdots+D_{j+p-1}$ be a chain of $p\ge 1$ rational curves contained in $C$ (where the indices are taken modulo $s$). Then $D'$ is of type {\bf a} or of type {\bf b} and  $D'$ is of type {\bf b} if $D'$ contains a subchain of type {\bf b}. In particular a cycle cannot contain two disjoint chains of type {\bf b}.
\end{Lem}
Proof: On induction on $p\ge 1$, the case $p=1$ being proved in Lemma \thesection.\ref{N2(2.5)}. Let $D'=D_j+\cdots+D_{j+p}$ be a chain of $p+1$ rational curves contained in $C$. By induction hypothesis $D''=D_j+\cdots+D_{j+p-1} $ is of type {\bf a} or {\bf b} and we have to consider $D'=D''+D_{j+p}$. Since 
$$(-2e_i-e_I)(-2e_j-e_J)\le 0\leqno{(\ast)}$$
not both are of type {\bf b}.
\begin{itemize}
\item If one, say $D''$ is of type {\bf b},  we have in $H_2(S,\bb Z)$,
$$D''= -2e_i-e_{I},\quad D_{j+p}= e_j-e_J$$
with $i\not\in I$, $j\not\in J$. Therefore, 
$$1=(-2e_i-e_I)(e_j-e_J)=-2e_ie_j- e_je_I + 2e_ie_J + e_Ie_J$$
and 
\begin{itemize}
\item either $j\in I$, $i\neq j$, $i\not\in J$, $I\cap J=\emptyset$, whence
$$D''+D_{j+p}= -2e_i-e_I+e_j-e_J=-2e_i - e_{I\cup J\setminus\{j\}}$$
is of type {\bf b},
\item or $j\not\in I$, $i=j$, $i\not\in J$ and $I\cap J$ contains one element, say $k$. Setting $I'=I\setminus\{k\}$ and $J'=J\setminus\{k\}$ we obtain
$$D''+D_{j+p}= -2e_k - e_{\{i\}\cup I'\cup J'}$$
is of type {\bf b}.
\end{itemize}
\item Both are of type {\bf a}
$$D'' = e_i-e_{I},\quad D_{j+p}= e_j-e_J$$
with $i\not\in I$, $j\not\in J$.  Clearly $i\neq j$ and
$$1=D''D_{j+p}=-e_ie_J - e_je_I + e_Ie_J,$$
then $i\in J$ or $j\in I$, say $j\in I$.\\
If $I\cap J=\emptyset$, $i\not\in J$, we set $I'=I\setminus\{j\}$ and $D''+D_{j+p}=e_i-e_{I'\cup J}$ is of type {\bf a}.\\
If $I\cap J=\{k\}$, $i\in J$ and we set $I'=I\setminus\{j,k\}$, $J'=J\setminus\{i,k\}$. Therefore $$D''+D_{j+p}=-2e_k-e_{I'\cup J'}$$
is of type {\bf b}.
\end{itemize}
In the following examples we use notations used in \cite{D1}, $a(S)=(a_i)_{i\in\bb Z}$ is the sequence of (opposite of) self-intersections of the rational curves in the universal covering space $\tilde S$ of $S$. This sequence is periodic of period $n$ and we overline a period. 

\begin{Ex}\label{Ex333} {\rm If $a(S)=(\overline{333})$ then we obtain a cycle $C=D_0+D_1+D_2$,
\begin{center}
\includegraphics[width=3cm]{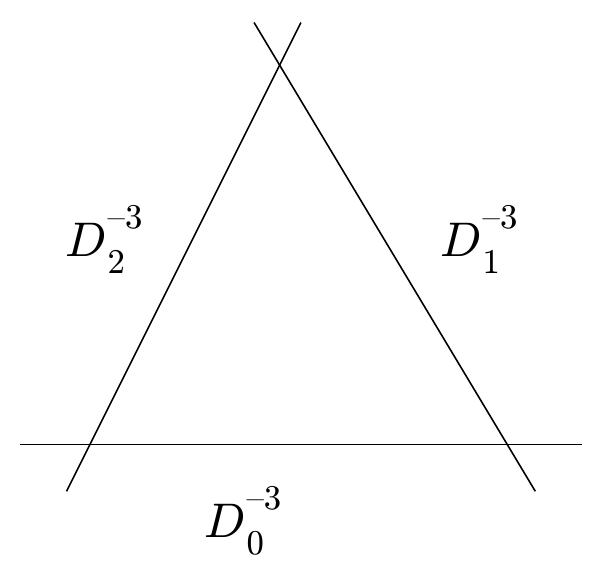}
\end{center}
All the curves are of type {\bf a}: $D_0= e_0-(e_1+e_2)$, $D_1= e_1-(e_2+e_0)$, $D_2= e_2-(e_0+e_1)$, $C=-(e_0+e_1+e_2)$,  however the chain $D_0+D_1= -2e_2$ is of type {\bf b} with $I$ empty.}
\end{Ex}
The previous example \thesection.\ref{Ex333} is generalized in the following lemma:
\begin{Lem}\label{internonvide} Let $D'=D_0+\cdots+D_{j-1}$ be a chain of $j\ge 2$ rational curves of type {\bf a}. The following two conditions are equivalent:\\
i) $D'$ is of type {\bf b} and all strict subchains are of type {\bf a}.\\
ii) $I_0\cap I_{j-1}=\{k\}$ with $k\not\in[0,j-1]$ and for any $0\le p<q\le j-1$, $(p,q)\neq (0,j-1)$ we have $I_p\cap I_q=\emptyset$.
\end{Lem}
Proof: Denote $D_i=e_i-e_{I_i}$.\\
If $j=2$, $D'=D_0+D_1= e_0-e_{I_0}+e_1-e_{I_1}$ is of type {\bf b}, then $I_0\cap I_1=\{k\}$ with $k\ge 2$. Conversely, if $I_0\cap I_1=\{k\}$ with $k\ge 2$, then $D_0+D_1$ is of type {\bf b}.\\
 Suppose that $j\ge 2$ and let $D'=D_0+\cdots+D_j=A+D_j$ be a chain of $j+1$ curves of type {\bf a} with $A=D_0+\cdots+D_{j-1}$.\\
If i) is satisfied, then since all subchains are of type {\bf a}, $A= e_0-e_{I_A}$, with $I_p\cap I_q=\emptyset$ when $0\le p<q\le j-1$. Moreover since $A+D_j$ is of type {\bf b}, then $I_A\cap I_j=\{k\}$, in particular $k\in I_j$ with $k>j$. Consider the subchain $D_1+\cdots+D_j$ which is of type {\bf a} as subchain. Then for $1\le p<q\le j$, $I_p\cap I_j=\emptyset$. Therefore $I_0\cap I_j=\{k\}$ with $k> j$.\\
Conversely if ii) is satisfied, subchains are of type {\bf a} by induction. If $I_0\cap I_j=\{k\}$ with $k> j$, then $D'=-2e_k-e_K$.\hfill 
$\Box$\\

Notice that the example \thesection.\ref{Ex333} is not a counterexample of the following lemma since the known example is an odd Inoue-Hirzebruch surfaces (= half-Inoue surface) with $a(S)=(\overline{s_1s_1s_1})=(\overline{333})$. In this case $[H_1(S,\bb Z):H_1(C,\bb Z)]=2$.
\begin{Lem} \label{lemme1.6} Let $C=D_0+\cdots+D_{s-1}$ be a cycle of $s\ge 2$ rational curves, with $D_i= e_i-e_{I_i}$, $i=0,\ldots,s-1$. If $H_1(C,\bb Z)=H_1(S,\bb Z)$, then
$$\forall\  i,j,\  i\neq j \Rightarrow  I_i\cap I_j=\emptyset.$$
\end{Lem}
Proof: Let $p:S'\to S$ be a 2-sheeted covering of $S$. By assumption $C'=p^\star(C)$ is a cycle $C'$ with $2s$ rational curves $$C'=D'_0+\cdots+D'_{s-1}+D'_n+\cdots+D'_{n+s-1}$$
 with $D'_0D'_1=\ldots=D'_{s-1}D'_n=\ldots=D'_{n+s-1}D'_0=1$. Suppose that $I_i\cap I_j=\{k\}$. Denoting $e'_0,\ldots,e'_{2n-1}$ the $2n$ Donaldson classes which trivialize the intersection form, we have $p^\star e_k = e'_k+e'_{n+k}$,
 $$p^\star(D_i)=D'_i+D'_{n+i}= (e'_i-e'_{I_i})+(e'_{n+i}-e'_{n+I_i}),$$
 $$p^\star(D_j)=D'_j+D'_{n+j} = (e'_j-e'_{I_j})+(e'_{n+j}-e'_{n+I_j}).$$ 
 If $I'_i\cap I'_j=\{k\}$ then  $I'_{n+i}\cap I'_{n+j}=\{n+k\}$. By Lemma \thesection.\ref{internonvide} we should obtain two disjoint chains of type {\bf b} and this is impossible by Lemma \thesection.\ref{chaineb}.

\hfill $\Box$
 
\begin{Lem} \label{Cyclea} Let $S$ be a minimal surface with a cycle $$C=D_0+\cdots+D_{s-1}= -(e_r+\cdots+e_{n-1})$$
 of $s\ge 2$ rational curves. Suppose that $H_1(C,\bb Z)=H_1(S,\bb Z)$ and
all curves are of type {\bf  a}. Then $s=r$, i.e.
$$\sharp(C)-C^2=b_2(S)$$
and $(I_i)_{0\le s-1}$ is a partition of $[0,n-1].$
\end{Lem}
Proof: 1) {\bf \boldmath If $s=2$},  $D_0 = e_{i_0}-e_{I_0}$, $D_1 = e_{i_1}-e_{I_1}$. We have
$$2=D_0D_1=-e_{i_0}e_{I_1} -e_{i_1}e_{I_0}+ e_{I_0}e_{I_1},$$
whence $i_1\in I_0$, $i_0\in I_1$, $I_0\cap I_1=\emptyset$. Setting $I'_0=I_0\setminus\{i_1\}$ and $I'_1=I_1\setminus\{i_0\}$, we obtain
$$D_0 = e_{i_0}-e_{i_1}-e_{I'_0}, \quad D_1 = e_{i_1}-e_{i_0}-e_{I'_1}, \quad {\rm with}\; I'_0\cap I'_1=\emptyset, \; \{i_0,i_1\}\cap (I'_0\cup I'_1)=\emptyset.$$
Therefore
$$-(e_r+\cdots+e_{n-1}) = C=D_0+D_1= -(e_{I'_0}+e_{I'_1})$$
i.e. $I'_0\cup I'_1=[r,n-1]$. Let $I=\{i_0,i_1\}\cup I'_0\cup I'_1$ and $I'=[0,n-1]\setminus I$. Suppose that $I'$ is non empty. Of course, $e_{I'}.D_0=e_{I'}.D_1=0$, whence by lemma \thesection.\ref{N2(2.4)} iii), $I'=[0,r-1]$ which is impossible. Therefore $I'=\emptyset$ and $I_0\cup I_1=\{i_0,i_1\}\cup I'_0\cup I'_1=[0,n-1]$, i.e. $r=2$ and 
$$\sharp(C)-C^2=2+(n-2)=b_2(S).$$
2) {\bf \boldmath If $s\ge 3$}, $D_j = e_{i_j}-e_{I_j}$, $j=0,\ldots,s-1$ and we may suppose that the numbering is such that
$$D_0D_1=\ldots=D_{s-2}D_{s-1}=D_{s-1}D_0=1.$$
By Lemma \thesection.\ref{lemme1.6}, $I_i\cap I_j=\emptyset$ if $i\neq j$.  Moreover the equality
$$1=D_jD_{j+1}=-e_{i_j}e_{I_{j+1}} - e_{i_{j+1}}e_{I_{j}}$$
implies that $i_j\in I_{j+1}$ or $i_{j+1}\in I_j$.\\
Changing if necessary the numbering (i.e. we number the curves following the other orientation of the cycle) we suppose that $i_1\in I_0$, then a straightforward proof by induction shows that $i_{j+1}\in I_j$ for $j\ge 0$.\\
 Setting $I'_j=I_j\setminus\{i_{j+1}\}$, we have
$$\begin{array}{lcl}
-(e_r+\cdots+e_{n-1})&=& C=D_0+\cdots +D_{s-1}= \displaystyle\sum_{j=0}^{s-1}(e_{i_j} - e_{I_j})\\ 
&&\\
&=& - (e_{I'_0}+\cdots+e_{I'_{s-1}})
\end{array}$$
therefore
$$[r,n-1]=I'_0 \cup \cdots \cup I'_{s-1} \quad {\rm with}\quad \{i_0,\ldots,i_{s-1}\}\cap\Bigl( I'_0 \cup \cdots \cup I'_{s-1}\Bigr)=\emptyset$$
Setting $I=\{i_0,\ldots,i_{s-1}\}\cup I'_0 \cup \cdots \cup I'_{s-1}$ and $I'=[0,n-1]\setminus I$, we have $L_{I'}D_j=0$ for all $j=0,\ldots,s-1$, then if $I'\neq\emptyset$, lemma \thesection.\ref{N2(2.4)} would imply once again that $I'=[0,r-1]$ and this is impossible, hence $I'=\emptyset$ and $\{i_0,\ldots,i_{s-1}\}\cup I'_0 \cup \cdots \cup I'_{s-1}=[0,n-1]$, i.e. $r=s$ and $\{i_0,\ldots,i_{r-1}\}=[0,r-1]$. Finally
$$\sharp(C)-C^2=s+(n-r)=n=b_2(S).$$
\hfill $\Box$\\

We determine now the surfaces for which there exists a rational curve of type {\bf b}. Following notations of \cite{D1}, $s_n$ is the sequence of $n$ integers $s_n=(n+2,2,\ldots,2)$.

\begin{Lem} \label{bsn} Let $S$ be a minimal surface with a cycle $C=D_0+\cdots+D_{s-1}= -(e_r+\cdots+e_{n-1})$ of $s\ge 2$ rational curves.   If there exists a rational curve of type {\bf b}, say $D_0$, then $s=n$, $r=0$, 
$$D_0 = -2e_1-e_{[2,n-1]},\quad D_1 = e_1-e_2,\ldots,D_{n-1} = e_{n-1}-e_0,$$
 $S$ is the odd Inoue-Hirzebruch surface such that $a(S)=(\overline{ s_n})$. In particular $[H_1(C,\bb Z):H_1(S,\bb Z)]=2$.\\
\begin{center}
\includegraphics[width=3cm]{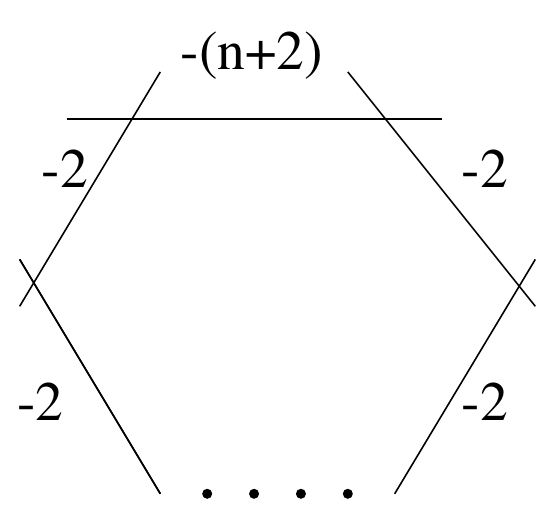}
\end{center}
\end{Lem}
Proof: First we prove  that $[H_1(C,\bb Z):H_1(S,\bb Z)]=2$ by contradiction. Suppose the contrary, then by \cite{N1}(2.13),  $[H_1(C,\bb Z):H_1(S,\bb Z)]=1$ and $C'=p^\star(C)$ is a cycle with $2s$ rational curves. Setting
$$p^\star(D_i)=D'_i+D'_{n+i}, \quad p^\star(e_i)=e'_i+e'_{n+i}$$
we have
$$D'_0D'_1=\cdots=D'_{s-2}D'_{s-1}=D'_{s-1}D'_n=D'_nD'_{n+1}=\cdots=D'_{n+s-1}D'_0=1.$$
If we choose the numbering such that $D_0 = -2e_{1}-e_{I_0}$, 
$$D'_0+D'_{n} = -2(e'_1+e'_{n+1})-e'_{I_0}-e'_{n+I_0}.$$ 
Since $D'_0$ and $D'_{n}$ are of the same type and since they cannot be both of type {\bf b} by lemma \thesection.\ref{chaineb}, both are of type {\bf a} and there exists an index $i_0$ such that
$$D'_0 = e'_0-e'_n-e'_1-e'_{n+1}-e'_{I_0}, \quad D'_n=e'_n- e'_0 - e'_1-e'_{n+1}-e'_{n+I_0}.$$
 By lemma \thesection.\ref{Cyclea}, 
$$\sharp(C')-C'^2=b_2(S')$$
Besides
$$\begin{array}{ccl}
{C'}^2 & = &\dps \sum_{i=0}^{s-1}({D'}_i^2 + {D'}_{n+i}^2) + 4s = 2\left(-4-\Card(I_0)+\sum_{i=1}^{s-1}(-1-\Card(I_i)) +2s \right)\\
&&\\
&=&2\dps\left( -3+s-\sum_{i=0}^{s-1}\Card(I_i)\right).
\end{array}$$
and by lemma \thesection.\ref{Cyclea},
$$[0,2n-1]=\{0,n,1,n+1\}\cup \bigcup_{i=0}^{s-1}(I_i\cup n+I_i)$$
therefore $[0,n-1]=\{0,1\}\cup \bigcup_{i=0}^{s-1}I_i$ and
$$n=2+\sum_{i=0}^{s-1}\Card(I_i).$$
These equalities yield
$$\begin{array}{ccl}
2b_2(S)&=&b_2(S')=\dps\sharp(C')-{C'}^2 = 2s+2\left( 3-s+\sum_{i=0}^{s-1}\Card(I_i)\right)\\
&&\\
&=&\dps 2\left( 1+2+\sum_{i=0}^{s-1}\Card(I_i)\right)=2(1+b_2(S)),
\end{array}$$
which is impossible.\\
Now, since $[H_1(C,\bb Z):H_1(S,\bb Z)]=2$,  a double covering $p:S'\to S$ yields a surface $S'$ with two cycles of rational curves, hence by \cite{N1}, $S'$ is an even Inoue-Hirzebruch surface and $S$ is an odd Inoue-Hirzebruch surface, in particular a Kato surface \cite[thm 3.8.]{D2}   By the explicit description of the self-intersection of the curves \cite{D2}, the only possible curve $D_0$ of type {\bf b} is $D_0 = e_0-e_{[1,n-1]}-e_0-e_1=-2e_1-e_{[2,n-1]}$, therefore $D_0^2=-(n+2)$ and $a(S)=(\overline{ s_n})$.\\\hfill$\Box$\\
\begin{Rem} In the previous computation subtraction $-e_k$ for a Kato surface means that we consider the strict transform by the blow-up by the exceptional curve $C_k$.
\end{Rem}

We notice that if $f:S'\to  S$ if a finite covering with $p$ sheets, then for $C'=f^\star(C)$,$$\sharp(C')=p\sharp(C),\quad C'^2=pC^2,\quad b_2(S')=pb_2(S).$$
Therefore the following result still holds for a rational curve with double point. We have proved:

\begin{Th}\label{EqTopCourbe} Let $S$ be a minimal surface with a cycle $C=D_0+\cdots+D_{s-1} = -(e_r+\cdots+e_{n-1})$ of $s\ge 1$ rational curves. \\
1) $S$ is an odd Inoue-Hirzebruch surface if and only if $[H_1(C,\bb Z):H_1(S,\bb Z)]=2$. In this case $s=n$, $r=0$ and 
$$\sharp(C)-C^2=2b_2(S).$$
2) If $H_1(C,\bb Z)=H_1(S,\bb Z)$, then $s=r$,
$$\sharp(C)-C^2=b_2(S),$$
all curves are of type ${\bf a}$ i.e. $D_i = e_i-e_{I_i}$, $i\not\in I_i$, and $(I_i)_{0\le i\le s-1}$ is a partition of $[0,n-1]$.

\end{Th}
We check these general results on Kato surfaces:
\begin{Ex}{\rm Let $S$ be the Inoue-Hirzebruch surface such that $$a(S)=(\overline{s_3s_1s_2})=(\overline{522342)}).$$
\begin{center}
\includegraphics[width=14cm]{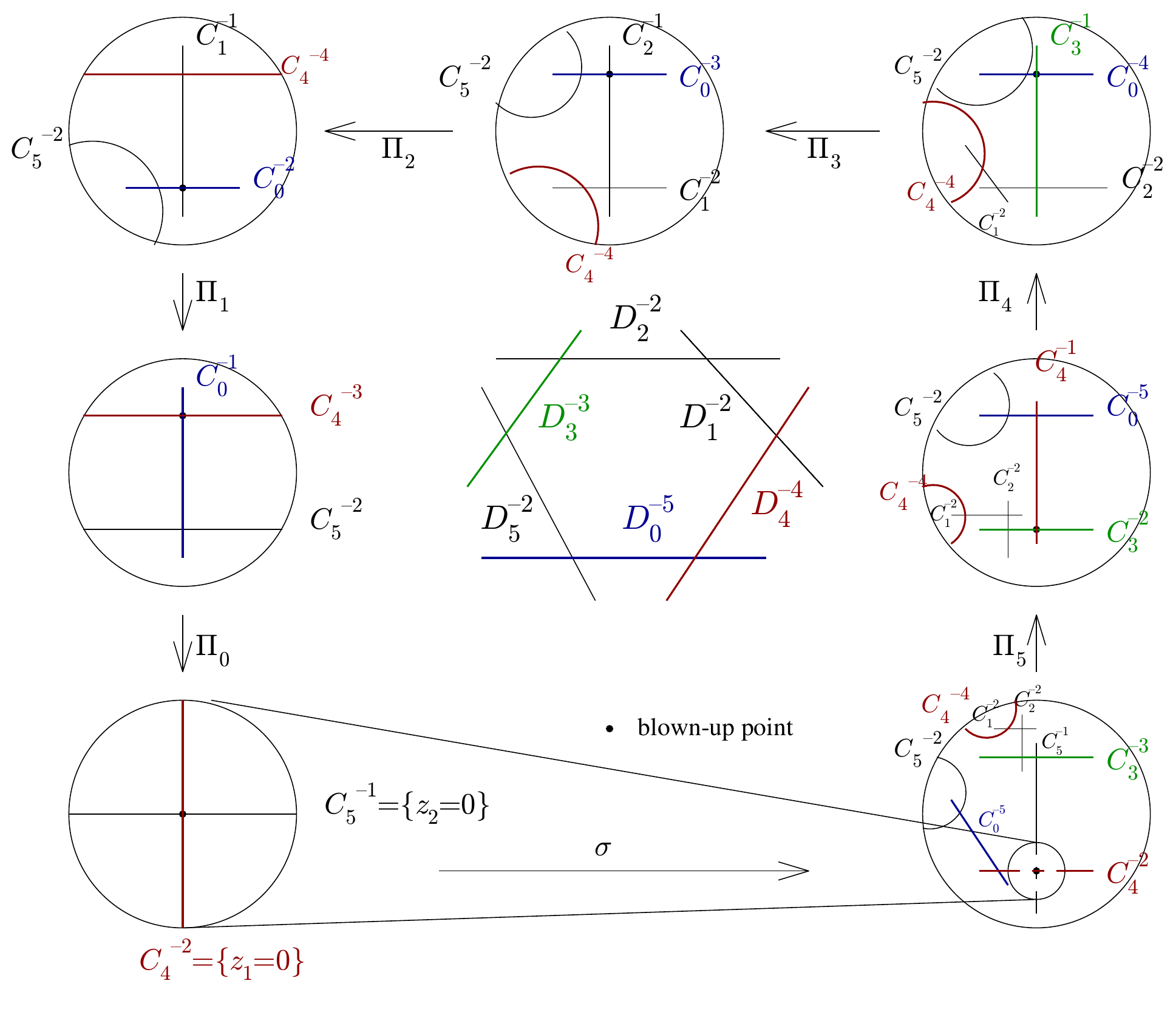}
\end{center}
We denote by $D_i$ the curve induced in $S$ by $C_i$. Then
$$\begin{array}{lcllcl}
D_0&=&e_0-(e_1+e_2+e_3+e_4)&D_1&= &e_1-e_2\\
D_2&=&e_2-e_3&D_3&=&e_3-(e_4+e_5)\\
D_4&=&e_4-(e_5+e_0+e_1)&D_5&=&e_5-e_0
\end{array}$$
Notice that there are two copies of each index in $\dps\coprod_{i=0}^5 I_i$.\\
Moving the blown-up point along $C_0$, we obtain a surface 
$$a(S)=(\overline{s_3s_1s_1r_1})=(\overline{522332)}),$$
with a regular sequence, hence with a branch (see \cite{D1})
\begin{center}
\includegraphics[width=14cm]{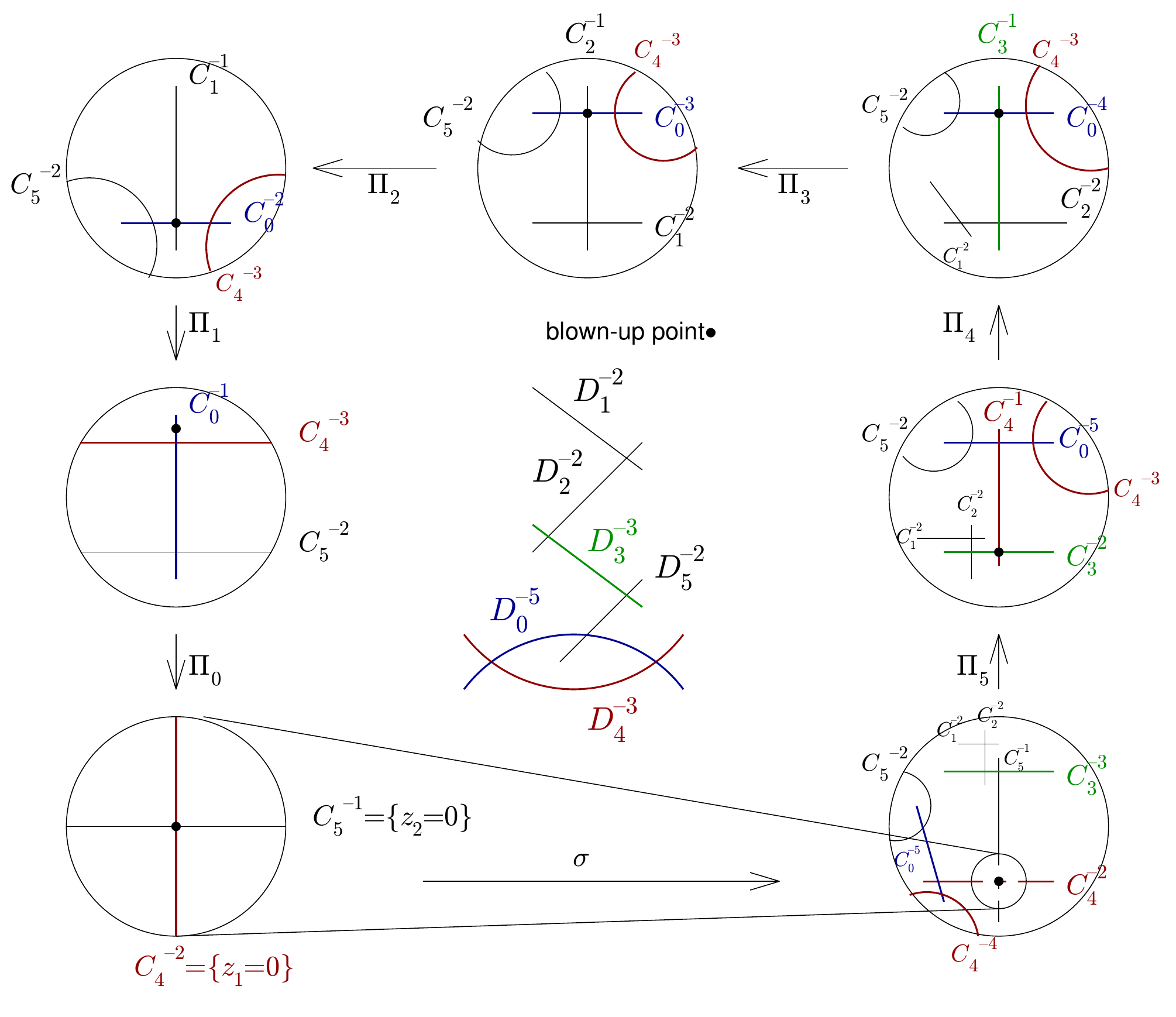}
\end{center}
$$\begin{array}{lcllcl}
D_0&=&e_0-(e_1+e_2+e_3+e_4)&D_1&= &e_1-e_2\\
D_2&=&e_2-e_3&D_3&=&e_3-(e_4+e_5)\\
D_4&=&e_4-(e_5+e_0)&D_5&=&e_5-e_0
\end{array}$$
The cycle is $C=D_0+D_4$,\  $\sharp(C)-C^2=6=b_2(S)$, $I_0\cap I_4=\emptyset$ and $I_0\cup I_4=[0,5]$.}
\end{Ex}

It is well-known that any cycle may be realized in a even Inoue-Hirzebruch surface $S'$, but in the following corollary we show that the second Betti number is the same:
 \begin{Cor} Let $S$ be a minimal surface with a cycle $C=D_0+\cdots+D_{r-1}\sim -(e_r+\cdots+e_{n-1})$ of $r\ge 1$ rational curves and let $n=b_2(S)$.  Then there is an even Inoue-Hirzebruch surface $S'$ with $b_2(S')=b_2(S)$ such that one of the two cycles, say $A=D'_0+\cdots+D'_{r-1}$ has the same intersection matrix.
\end{Cor}
Proof:  We have by the previous theorem 
$$b_2(S')=-C'^2+\sharp(C')=-C^2+\sharp(C)=b_2(S).$$
\hfill$\Box$

\begin{Cor} \label{classescourbescycle} Let $S$ be a surface with a cycle $C$ of rational curves. When the curve numbered $0$ is  fixed, there is a numbering of the rational curves of the cycle and of the homological classes $e_i$, $i=0,\ldots, n-1$ such that (with $\a_0=0$) 
$$C=D_0+\cdots+D_{\a_i}+\cdots+D_{\a_{s-1}}$$

and in $H_2(S,\bb Z)$, for $i=0,\ldots,s-1$,
$$\begin{array}{l}
D_{\a_0}=D_0= e_{0}-(e_1+\cdots+e_{\a_1}),\hfill -D_0^2=\a_1+1,\\
\\
 D_{\a_1}= e_{\a_1}-(e_{\a_1+1}+ \cdots + e_{\a_2}), \hfill -D_{\a_1}^2=\a_2-\a_1+1,\\
\hfill\vdots\hfill\hfill\vdots\hfill\hfill\\
 D_{\a_{i}}= e_{\a_{i}}-(e_{\a_{i}+1}+\cdots+e_{\a_{i+1}}), \hfill -D_{\a_{i}}^2=\a_{i+1}-\a_i+1\\
\hfill \vdots\hfill\hfill\vdots\hfill\hfill\\
 D_{\a_{s-1}}= e_{\a_{s-1}}-(e_{\a_{s-1}+1}+\cdots+e_{n-1}+e_{0}), \quad\hfill -D_{\a_{s-1}}^2=n-\a_{s-1}+1.
 \end{array}$$
 Notice that 
 $$D_0D_{\a_1}=\cdots=D_{\a_{i}}D_{\a_{i+1}}=\cdots=D_{\a_{s-1}}D_0=1.$$
 \end{Cor}

In other words, denote by $\mathfrak G\subset [0,n-1]$ the subset of indices of rational curves in $C$, i.e.
$$\mathfrak G=\{\a_0,\a_1,\ldots,\a_{s-1}\}$$
then for $I=I_C=[0,n-1]\setminus\mathfrak G$, $C=-e_{I}\in H_2(S,\bb Z)$. If there is a curve of type {\bf b}, $S$ is a Inoue-Hirzebruch surface by lemma \thesection.\ref{bsn}, in particular $S$ is a Kato surface, therefore in the sequel we shall suppose that all curves are of type {\bf a} and we shall write any curve of $C$ as
$$D_{\a_{i}}= e_{\a_{i}}-e_{I'_{\a_{i}}}  -e_{\a_{i+1}}$$
with $I'_{\a_{i}}=I_{\a_{i}}\setminus\{\a_{i+1}\}$, then
$$I_C=\bigcup_{i=0}^{s-1} I'_{\a_0+\cdots+\a_i}.$$

\section{(Co)homological class of the maximal divisor}

The following theorem shows that the class of the maximal divisor of a surface with a cycle of rational curves is of the same type as the class of a cycle. Moreover the maximal divisor looks like the one of a Kato surface: if there is one connected component there is a cycle and trees attached to different curves.
\begin{Lem} \label{ejectioncourbeexc} Let $S$ be a compact complex surface in class VII$_0^+$ endowed with a cycle  $C=D_{\a_0}+\cdots+D_{\a_i}+\cdots+D_{\a_{s-1}}$ of $s$ rational curves, $1\le s\le n$, of class $-e_{I_C}$, $I_C\subset [0,n-1]$. Let $D_{\a_i}=e_{\a_i}-(e_{\a_1+1}+\cdots+e_{\a_{i+1}})$ be the class of $D_{\a_i}$, i.e. $I_C=[0,n-1]\setminus\{\a_0,\ldots,\a_{s-1}\}$. Then there exists a deformation over the disc $\cal S\to \D$ which smoothes the point $D_{\a_i}\cap D_{\a_{i+1}}$, i.e. there is a flat families $\cal C\to\D$ such that 
\begin{itemize}
\item $S_0\simeq S$, $S_t$ is not minimal, contains exactly one exceptional curve of the first kind of class $e_{\a_{i+1}}$,
\item $C_0\simeq C$,  $C_t=D_{\a_0}+\cdots+D'_{\a_i}+D_{\a_{i+2}}+\cdots+D_{\a_{s-1}}$, $t\neq 0$, is of class $-e_{I_C}$ and $D'_{\a_i}=D_{\a_i}+D_{\a_{i+1}}=e_{\a_i}-(e_{\a_i+1}+\cdots+e_{\a_{i+1}-1}+e_{\a_{i+1}+1}+\cdots+e_{\a_{i+2}})$.
\end{itemize}
\end{Lem}
Proof: By Karras \cite{Karras77}, it is possible to smooth the singular point $D_{\a_i}\cap D_{\a_{i+1}}$ in order to obtain a cycle $C_t$ of rational curves (an elliptic curve if $s=1$) of the same class, in particular of the same self-intersection, but with one fewer rational curve.  This local deformation can be globally realized by a deformation of surfaces. In fact, we have the exact sequence of sheaves
$$0\to \T_S(-\log C)\to \T_S\to J_C\to 0$$
with $J_C=\T_S/\T_S(-\log C)$. Since for any curve $D_{\a_i}$ of $C$, $D_{\a_i}^2\le -2$, $H^0(S,J_C)=0$. Besides with \cite[lemma (4.3)]{N1} and \cite[thm (1.3)]{N2}, the long exact sequence of cohomology gives the exact sequence
$$0\to H^1(S,\T_S(-\log C))\to H^1(S,\T_S)\to H^1(U,\T_U)\to 0$$
for a strictly pseudoconvex neighbourhood of $C$.\\
We have in $H_2(S,\bb Z)$
$$D'_{\a_i}=D_{\a_i}+D_{\a_{i+1}}$$
By the formula of theorem \ref{EqTopCourbe}, there is exactly one exceptional curve $E_t$ of the first kind of class the missing one in the expressions of classes of the curves in $C_t$, i.e. $E_t=e_{\a_{i+1}}$.\hfill$\Box$

\begin{Th} \label{compconnexe} Let $S$ be a compact complex surface in class VII$_0^+$ endowed with exactly one cycle  $C=D_0+\cdots+D_{\a_i}+\cdots+D_{\a_{s-1}}$ of $s$ rational curves, $1\le s<n$, of class $-e_{I_C}$. Let $D=C+A$ be the maximal connected reduced divisor containing $C$ with $A\neq 0$. Let $D_{\a_s},\ldots,D_{\a_{s+q-1}}$,  be the irreducible components of $A$ numbered so that $(C\cup D_{\a_s}\cup\cdots\cup D_{\a_{s+k-1}})\cap D_{\a_{s+k}}\neq \emptyset$, $k=0,\ldots,q-1$ (i.e. the unions remain connected).  Then:
\begin{enumerate}
\item There is a family of subsets $\mathfrak D_k\subset [0,n-1]$, $k\ge 0$,  such that $\mathfrak D_0=I_C$ and
$$D=C+A=-e_{\mathfrak D_{q-1}}\in H_2(S,\bb Z),$$
\item Each connected component $A_i$ of $A$, called a tree, is a chain of rational curves,
\item Two different trees $A_j$, $A_k$, meet different curves of the cycle.
\end{enumerate}
\end{Th}
Proof:  Since $p<n$, $S$ is not an odd Inoue-Hirzebruch surface and each class of curve has the form $e_i-e_{I_i}$. If $D_{\a_{s+k}}=e_{\a_{s+k}}-e_{I_k}$, we have for $k\ge 0$,
$$\begin{array}{lcl}
1\le D_{\a_{s+k}}.(C+D_{\a_s}+\cdots+D_{\a_{s+k-1}})&=&-(e_{\a_{s+k}}-e_{I_k}).e_{\mathfrak D_{k-1}}\\
&&\\
&=&-e_{\a_{s+k}}.e_{\mathfrak D_{k-1}} + e_{I_{k}}.e_{\mathfrak D_{k-1}}\le  1
\end{array}$$
therefore 
\begin{enumerate}
\item $\a_{s+k}\in \mathfrak D_{k-1}$,
\item  $I_{k}\cap \mathfrak D_{k-1}=\emptyset$
\item there exists  a unique index $i\in \{0,\ldots, s+k-1\}$ such that $D_{\a_{s+k}}.D_{\a_i}=1$  i.e.  $D_{\a_{s+k}}$ meets only one curve at one point.
\end{enumerate}
 Finally
$$C+D_{\a_s}+\cdots+D_{\a_{s+k}}=-e_{\mathfrak D_{k-1}} + (e_{\a_{s+k}}-e_{I_{k}})= -e_{\mathfrak D_{k}}$$
where $\mathfrak D_{k}=(\mathfrak D_{k-1}\setminus \{\a_{s+k}\})\cup I_{k}.$\\

 At each step $k\ge 0$, the indices of curves are in
$$\{\a_0,\ldots,\a_{s-1},\a_s,\ldots,\a_{s+k-1}\}$$
and $D_{\a_{s+k}}$ meets a curve $D_{\a_i}$ of the cycle without any other tree or the top of a tree. In fact, suppose the contrary:
\begin{itemize}
\item If $D_{\a_{s+k_1}}=e_{\a_{s+k_1}}-e_{I_{s+k_1}}$ and $D_{\a_{s+k_2}}$ meet the same curve in the cycle, say $D_{\a_0}$. Then we have for $D_{\a_{s+k_1}}$ the system
$$\begin{array}{cclcl}
1&=&D_{\a_{0}}.D_{\a_{s+k_1}}&=&\lbrack e_{\a_{0}}-(e_{\a_{0}+1}+\cdots+e_{\a_{1}}\rbrack.\lbrack e_{\a_{s+k_1}}-e_{I_{s+k_1}}\rbrack\\
0&=&D_{\a_1}.D_{\a_{s+k_1}}&=&\lbrack e_{\a_1}-(e_{\a_1+1}+\cdots+e_{\a_{2}})\rbrack.\lbrack e_{\a_{s+k_1}}-e_{I_{s+k_1}})\rbrack\\
&&\ldots&&\ldots\\
0&=&D_{\a_{j}}.D_{\a_{s+k_1}}&=&\lbrack e_{\a_{j}}-(e_{\a_{j}+1}+\cdots+e_{\a_{j+1}})\rbrack.\lbrack e_{\a_{s+k_1}}-e_{I_{s+k_1}})\rbrack\\
&&\ldots&&\ldots\\
0&=&D_{\a_{s-1}}.D_{\a_{s+k_1}}&=&\lbrack e_{\a_{s-1}}-(e_{\a_{s-1}+1}+\cdots+e_{\a_{0}})\rbrack.\lbrack e_{\a_{s+k_1}}-e_{I_{s+k_1}})\rbrack
\end{array}$$
Recall that if $D'=e_{i'}-e_{I'}$ and $D''=e_{i''}-e_{I''}$ are distinct curves, then $i'\neq i''$.

\begin{enumerate}
\item If $\a_{s+k_1}\in\{\a_0+1,\ldots,\a_1\}$ the system is equivalent to 
$$\begin{array}{ccl}
0&=& -e_{\a_0}.e_{I_{s+k_1}}+(e_{\a_0+1}+\cdots+e_{\a_{1}}).e_{I_{s+k_1}}\\
0&=& -e_{\a_1}.e_{I_{s+k_1}} +(e_{\a_1+1}+\cdots+e_{\a_{2}}) e_{I_{s+k_1}}\\
&&\ldots\\
0&=&-e_{\a_{s-1}}.e_{I_{s+k_1}} +(e_{\a_{s-1}+1}+\cdots+e_{\a_{0}}) e_{I_{s+k_1}}
\end{array}$$
We have $I_{s+k_1}\subset\{\a_0,\ldots,\a_{s-1}\}$ (recall condition $(2)$), then there are two cases
\begin{itemize}
\item if $\a_0\not\in I_{s+k_1}$ then $\a_1\not\in I_{s+k_1}$,\ldots, $\a_{s-1}\not\in I_{s+k_1}$ and $I_{s+k_1}=\emptyset$ which is impossible. Contrarily 
\item if $\a_0\in I_{s+k_1}$, then the equations imply $\a_1\in I_{s+k_1}$, \ldots, $\a_{s-1}\in I_{s+k_1}$ and 
$$I_{s+k_1}=\{\a_0,\ldots,\a_{s-1}\},\quad \Card I_{s+k_1}=s\ge 1.$$
\end{itemize}
\item If $\a_{s+k_1}\in\{\a_j +1,\ldots,\a_{j+1}\}$, $j>0$, the system is equivalent to 
$$\begin{array}{ccl}
1&=& -e_{\a_0}.e_{I_{k_1}}+(e_{\a_0+1}+\cdots+e_{\a_{1}}).e_{I_{s+k_1}}\\
0&=& -e_{\a_1}.e_{I_{k_1}} +(e_{\a_1+1}+\cdots+e_{\a_{2}}). e_{I_{s+k_1}}\\
&&\ldots\\
0&=& -e_{\a_j}.e_{I_{k_1}}+1+(e_{\a_j+1}+\cdots+e_{\a_{j+1}}).e_{I_{s+k_1}}\\
0&=& -e_{\a_{j+1}}.e_{I_{k_1}}+(e_{\a_{j+1}+1}+\cdots+e_{\a_{j+2}}).e_{I_{s+k_1}}\\
&&\ldots\\
0&=&-e_{\a_{s-1}}.e_{I_{k_1}} +(e_{\a_{s-1}+1}+\cdots+e_{\a_{0}}) e_{I_{s+k_1}}
\end{array}$$
and since $I_{s+k_1}\subset\{\a_0,\ldots,\a_{s-1}\}$
$$I_{s+k_1}=\{\a_{j+1},\ldots,\a_{s-1},\a_0\},\quad \Card I_{s+k_1}=s-j\ge 1.$$
\end{enumerate}
We have the same conditions for $D_{\a_{s+k_2}}$ therefore an easy calculation shows that $D_{\a_{s+k_1}}.D_{\a_{s+k_2}}<0$ which is impossible.
\item We have to check now that each connected component $A_i$ of $A$ is a chain of rational curves. We prove that by induction on the number of the rational curves of the cycle.\\
{\bf If \boldmath$p=1$}, we can number the Donaldson classes such that 
$$D_0=e_0-(e_1+\cdots+e_{n-1}+e_0)= -(e_1+\cdots+e_{n-1}).$$
 If a curve $D_1=e_k-e_{I_1}$ meets $D_0$,
$$\begin{array}{lll}
1\le D_0.D_1&=&-(e_1+\cdots+e_{n-1}).(e_k-e_{I_1})\\
&=&-(e_1+\cdots+e_{n-1}).e_k + (e_1+\cdots+e_{n-1}).e_{I_1}\le 1
\end{array}$$
therefore $1\le k\le n-1$ and $I_1 \cap\{1,\ldots,n-1\}=\emptyset$. It means that $I_1=\{0\}$. If another curve $D'=e_{k'}-e_{I'}$ meets $D_0$ we have also $I'=\{0\}$ and $D_1.D'=(e_k-e_0).(e_{k'}-e_0)=-1$ which is impossible. We choose the numbering such that $k=1$ and $D_1=e_1-e_0$. We have $D_0+D_1=-(e_0+e_2+\cdots+e_{n-1})$. Let $D_2=e_l-e_{I_2}$ with a non-empty intersection with $D_0+D_1$. We have
$$1\le (D_0+D_1).D_2= -(e_0+e_2+\cdots+e_{n-1}).(e_l-e_{I_2})\le 1.$$
Hence $l\in \{0,2,\ldots,n-1\}$ and $I_{2}=\{1\}$, i.e. if numbering is choosen such that $l=2$, $D_2=e_2-e_1$. By induction we obtain a chain.\\
{\bf If \boldmath $p>1$} Suppose that we have the result for cycles of $p-1\ge 1$ curves, and let $C=D_{\a_0}+\cdots+D_{\a_p-1}$ be a cycle of $p$ curves. By lemma \ref{ejectioncourbeexc}, there is a deformation which smoothes a singular point of the cycle, say $D_{\a_0}\cap D_{\a_1}$, then self-intersection does not change and the number of curves decreases by one, therefore the deformed surface in not minimal (recall the formula $b_2(S)=-C^2+\sharp(C)$ for minimal surfaces) with exactly one exceptional curve of the first kind of class $e_{\a_1}$. If the new rational curve $D'_0$, of the same class as  $D_{\a_0}+ D_{\a_1}$, is intersected by one connected component, all connected components are chains by the induction hypothesis. If two connected components meet $D'_0$, then by the induction hypothesis one of the two is an exceptional divisor. Since there is only one exceptional curve of the first kind this component is also a chain.
\end{itemize}
\hfill$\Box$
\begin{Rem} It is difficult to detect more curves than those in the cycle, since the expressions of $C$ and of the maximal divisor $D=C+A$ in $H_2(S,\bb Z)$
$$C=-e_{[0,n-1]\setminus\mathfrak G}, \quad D=\G+A=-e_{\mathfrak D}$$
have the same type, i.e. the coefficients are all equal to $0$ or $-1$;  they have also  same genus $p_a(C)=p_a(C+A)=1$, where the arithmetic genus $p_a$ of a divisor $\D$ is defined by $p_a(\D)=1+\frac{1}{2}(K\D+\D^2)$.\\
For example, if $b_2=2$ the class of a cycle $C$ does not determine the maximal divisor: it may be either a cycle $C=-e_1$ with one curve and another cycle $C'=-e_0$ or a cycle $C=-e_1$ with a tree $A=e_1-e_0$.
\end{Rem}

\section{Connected components of the maximal divisor}

\begin{Lem} \label{classediviseursc} Let $S$ be a surface in class VII$_0^+$ containing a cycle $C$ of rational curves with $n=b_2(S)\ge 1$. We suppose that $S$ is not an odd Inoue-Hirzebruch surface. Let $H$ be a reduced divisor whose support is connected and simply connected. Then there is an integer $k\in\{0,\ldots,n-1\}$ and a subset $K\subset \{0,\ldots,n-1\}$, with $k\not\in K$, such that in $H_2(S,\bb Z)$
$$H=e_k-e_K.$$
\end{Lem}
Proof: Since the surface is not an odd Inoue-Hirzebruch surface the class of any rational curve is of the type $e_i-e_I$. We prove the result by induction on the number of irreducible components with the same arguments already used in theorem \ref{compconnexe}.\hfill$\Box$

\begin{Th}\label{2composantes} Let $S$ be a surface in class VII$_0^+$ with $b_2(S)=n$, containing a cycle $C=D_{\a_0}+\cdots+D_{\a_{s-1}}$ of $s$ rational curves and $D=C+A$ the connected component of the cycle. If there is another connected component  of curves $C'\neq 0$, then $C'$ is a cycle of rational curves, $A=0$ and $S$ is a Inoue-Hirzebruch surface.
\end{Th}
Proof: Let $I\subset \{0,\ldots,n-1\}$ such that $C=-e_I$. By Corollary \ref{classescourbescycle}, indices of the curves of $C$ belong to $\{\a_0,\ldots,\a_{s-1}\}$ and $I=[0,n-1]\setminus\{\a_0,\ldots,\a_{s-1}\}$. Suppose that $C'$ does not contain a cycle of rational curves. Then $C'$ is simply connected. Since $C'.D=0$,
Lemma \ref{classediviseursc}, shows that $C'=e_k-e_K$, with $k\in I$, $k\not\in K$, and $K\cap I=\{l\}$. Similarly if $A_i$ is a tree, since $A_i.C=1$, $A_i=e_{j_i}-e_{J_i}$ with $j_i\in I$ and $J_i\cap I=\emptyset$. Moreover since two branches $A_i=e_{j_i}-e_{J_i}$ and $A_{i'}=e_{j_{i'}}-e_{J_{i'}}$ do not meet $J_i\cap J_{i'}=\emptyset$.
We apply now theorem (1.4) of \cite{N2} with $H=A+C'$ in order to obtain a deformation $\cal S\to \D^s$ with flat families $\cal C$ and $\cal H$ where $H_t=H$, and $\cal C\to\D^s$ is the versal deformation of $C$, in particular generically $C_t$ is an elliptic curve of class $-e_I$ and $S_t$ contains exceptional divisors of the first kind $E_{\a_0}=e_{\a_0}$, \ldots, $E_{\a_{s-1}}=e_{\a_{s-1}}$. The only minimal surfaces in class VII$_0^+$ with an elliptic curve are Hopf surfaces and Inoue surfaces. Therefore, generically, $S_t$ is a blown-up Inoue surface or a blown-up Hopf surface with minimal model $S'_t$ of second Betti number $0\le b_2(S'_t)\le n-s$. There is a disc $\D$ over which $C_t$ is an elliptic curve, blown-up by exceptional divisors of the first kind $E_i=e_i$, $i\in I'$, $I'\subset I$, exceptional divisors $E_{\a_j}$, $0\le j\le s-1$ and $C'=e_k-e_K$. If $I'\neq I$, $S_t$ is a blown-up Inoue surface and contains also a cycle 
$$\G_t=\sum_{i\in I\setminus I'}D_i,$$
 with $D_i=e_i-e_{i+1}$ (when $I\setminus I'$ is properly numbered modulo $\Card(I)-\Card(I')$).
Since $k\in I$ there are two cases:
\begin{itemize}
\item $k\in I'$, then $C'.E_k=(e_k-e_K).e_k=-1$ impossible,
\item $k\in I\setminus I'$ then $C'.D_k=(e_k-e_K).(e_k-e_{k+1})=-1+e_K.e_{k+1}<0$ impossible
\end{itemize}
If $I'=I$ we are in the first case.\\
In any case we obtain a contradiction.\\
Now, $C'$ contains a cycle of rational curves, we have two cycles, we conclude by a theorem of Nakamura \cite[Thm 8.1]{N1} that $S$ is a Inoue-Hirzebruch surface, $C'$ is a cycle and $A=0$.

\hfill$\Box$

\section{Application: Surfaces with twisted logarithmic forms}
\begin{Lem} \label{cycle} Let $S$ be a surface in class VII$_0^+$ which is not a Enoki surface. If there exists a non-trivial twisted logarithmic $1$-form with polar set $D$, 
$$0\neq \o\in H^0(S,\O^1(\log D)\ot \cal L)$$
 with $\cal L\in H^1(S,\bb C^\star)$, then $D\neq 0$ and each connected component of $D$ contains a cycle of rational curves. In particular, there is at most two connected components and in this case $S$ is a Inoue-Hirzebruch surface.
\end{Lem}
Proof: By \cite[Thm 5.6]{AD18} there exists  a line bundle $\cal L\in Pic^0(S)$ such that $H^1(S,\O^1\ot \cal L)\neq 0$ only if $S$ is a Enoki surface, therefore $D\neq 0$. Let $D'$ be any connected component of the polar set $D$ of $\o$, then the intersection matrix is negative definite by \cite{E81}.  Let a spc neighbourhood $V'$ of $D'$ sufficiently small to be contractible onto a Stein space $V$ with an isolated singular normal point $x\in V$. Let $\pi:V'\to V$ be this contraction, $D'=\pi^{-1}(x)$ be the exceptional divisor. Let $U=V\setminus\{x\}$.  If $D'$ does not contain a cycle, $D'$ is simply connected then the restriction of $\cal L$ to $V'$ is trivial . By the theorem of Steenbrinck-Van Straten \cite[Cor. 1.4]{SVS85} the mapping induced by the differentiation
$$d:H^0(U,\O^1_U)/H^0(V',\O^1_{V'})\to H^0(U,\O^2_U)/H^0(V',\O^2_{V'}(\log D'))$$
is injective, which gives a contradiction if there is no cycle. If there is a cycle and another connected component of rational curves then there are two cycles by \ref{2composantes}.\hfill$\Box$

\begin{Prop} Let $S$ be a surface in class VII$_0^+$. If there exists a non-trivial twisted logarithmic $1$-form with polar set $D$, 
$$0\neq \o\in H^0(S,\O^1(\log D)\ot \cal L)$$
 with $\cal L\in H^1(S,\bb C^\star)$, then the vanishing divisor of $\o$ is trivial.
 \end{Prop}
 Proof: If $S$ is a Enoki surface then any twisted logarithmic or holomorphic $1$-form does not vanish. If $S$ is not a Enoki surface and there is a divisor $A$ where $\o$ vanishes, $A\cap D=\emptyset$ and there are two connected components of curves, therefore $S$ would be a Inoue-Hirzebruch surface. However on these surfaces $\o$ does not vanishes at all.\hfill$\Box$

\begin{flushright} 
Georges Dloussky\\ 
Aix Marseille Université, CNRS, Centrale Marseille, I2M, UMR 7373, 13453, Marseille, France\\
georges.dloussky@univ-amu.fr\\
\end{flushright}

\begin{thebibliography}{AUT}
\bibitem{AD18}  {\sc Apostolov V., \& Dloussky G.} On the Lee classes of locally conformally symplectic complex surfaces. J. of symplectic geometry 16 n$^0$4, pp 931-958 (2018).

\bibitem{D1} {\sc G. Dloussky }: Structure des surfaces de
Kato. {\em M\'emoire de  la S.M.F 112.$\rm n^{\circ}14$ (1984).} 
\bibitem{D2} {\sc G. Dloussky }: Une construction \'el\'ementaire des
surfaces d'Inoue-Hirzebruch. {\em Math. Ann. 280, (1988), 663-682.} 

\bibitem{D06} {\sc Dloussky G.} On surfaces of class VII$_0^+$ with numerically anticanonical divisor. {\em Am. J. Math. 128, 639-670 (2006)}

\bibitem{DOT3} {\sc G. Dloussky, K. Oeljeklaus, M. Toma }: Class 
VII$_{0}$ 
surfaces with $b_{2}$ curves. {\em Tohoku Math. J. 55 (2003), 283-309}

\bibitem{DON}{\sc S.K. Donaldson }: The orientation of Yang-Mills moduli space and 4-manifolds topology. {\em J. Differential Geometry 26 (1987) 397-428.}


\bibitem{E81} {\sc I. Enoki }: Surfaces of class VII$_0$ with curves. {\em T\^ohoku Math. J. 33, (1981), 453-492.} 

 

 
\bibitem{Karras77}{\sc Karras Ulrich} Deformations of cusp singularities, {\em Proc. of Symp. in Pure Math. 30, 1977, p37-43.}


\bibitem{N1} {\sc I. Nakamura}: On surfaces of class $\rm VII_0$ with
curves. {\em Invent. Math. 78,(1984), 393-443.} 
\bibitem{N2}{\sc I. Nakamura} On surfaces of class $VII_{0}$ with 
curves 
II.{\em  Tohoku Math. J. 42 (1990),  475-516.}

\bibitem{SVS85}{\sc Joseph Steenbrink \& Duco van Straten} Extendability of holomorphic differential forms near isolated hypersurface singularities {\em Abh. Math. Sem. Univ. Hamburg 55, 97-110 (1985)}


\bibitem{Te05}{\sc Teleman, A.}:  Donaldson theory on non-K\"ahlerian
surfaces and class VII surfaces with $b_2=1$, {\em Invent. math. 162, 493-521 (2005)}


\bibitem{Te09} {\sc A.~Teleman,}  Gauge theoretical methods in the classification of non-Kählerian surfaces, {\em Algebraic topology Ð old and new, Banach Center Publ . 85 (2009) 109Ð120.}

\bibitem{Te10} {\sc Teleman, A.}: Instantons and curves on class VII surfaces,  {\em  Annals of Math. (2) 172 (2010), n¡3, 1749Ð1804.}



\end{thebibliography}
\end{document}